\documentclass[11pt,a4paper,reqno]{amsart}

\usepackage{times} 
\usepackage{amsmath} 
\usepackage{amssymb}  
\usepackage{amsthm}
\usepackage{latexsym}
\usepackage{amsfonts,bbm}
\usepackage{xcolor}
\usepackage{mathtools}
\usepackage{enumerate}
\usepackage{cite}
\usepackage{tikz}
\usepackage{graphicx,caption}
\usepackage{todonotes}
\usepackage{float}
\usepackage{nicefrac}

\newtheorem{theorem}{Theorem} 
\newtheorem{assumption}[theorem]{Assumption}
\newtheorem{remark}[theorem]{Remark}
\newtheorem{lemma}[theorem]{Lemma} 
\newtheorem{example}[theorem]{Example} 
\newtheorem{definition}[theorem]{Definition} 

\usepackage[foot]{amsaddr}
 
\numberwithin{equation}{section}
\allowdisplaybreaks

\newcommand{\changed}[1]{\textcolor{black}{#1}}

\title{%
	Perturbations in PDE-constrained optimal control decay exponentially in space
}
\author{Simone Göttlich}\address{Chair of Scientific Computing,  School of Business Informatics and Mathematics, University of Mannheim, Germany.\\ Mail: \textsc{goettlich@uni-mannheim.de}}
\author{Manuel Schaller}\address{Faculty of Mathematics, Chemnitz University of Technology,\\ Mail: \textsc{manuel.schaller@math.tu-chemnitz.de}}
\author{Karl Worthmann}\address{Optimization-based Control Group, Institute of Mathematics, Technische Universität Ilmenau, Germany.\\ Mail: \textsc{karl.worthmann@tu-ilmenau.de}}

\thanks{K.W. gratefully acknowledges funding by the German Research Foundation (DFG) – Project-ID 507037103.}
\begin{document}
\maketitle
\begin{abstract}
For linear-quadratic optimal control problems (OCPs) governed by elliptic and parabolic partial differential equations (PDEs), we investigate the impact of perturbations on optimal solutions.
Local perturbations may occur, e.g., due to discretization of the optimality system or {disturbed} problem data. 
Whereas these perturbations may exhibit global effects in the uncontrolled case, we prove that the ramifications 
are exponentially damped in space under stabilizability and detectability conditions. 
To this end, we prove a bound on the optimality condition's solution operator that is uniform in the domain size. Then, this uniformity is used in a scaling argument to show the exponential decay of perturbations in space. 
We numerically validate and illustrate 
our results by solving OCPs involving 
Helmholtz, Poisson, and 
advection-diffusion-reaction equations. 

\end{abstract}


\smallskip
\noindent \textbf{Mathematics subject classications (2020).} 35Q93, 49K40, 93D23

\bigskip

\section{Introduction}\label{sec:introduction}
\noindent We analyze control of partial differential equations (PDEs) on bounded domains $\Omega\subset \mathbb {R}^d$, $d\in \mathbb{N}$, with Lipschitz boundary $\partial \Omega$. 
We consider elliptic linear quadratic optimal problems of the form
\begin{align}\label{eq:OCPintro}
\min_{(y,u)\in V\times U}\quad \frac12 \left(\|C(y-y_d)\|_Y^2 \right.&\left.+ \|u\|_U^2\right) \qquad\text{subject to}\qquad Ay + Bu = f,
\end{align}
where $V \subset H^1(\Omega)$ depends on the choice of
boundary conditions, $A\in L(V,V^*)$ 
is a surjective (not necessarily boundedly invertible) second-order differential operator in weak form and $B\in L(U,V^*)$, $C\in L(V,Y)$ are input and output operators with input space $U$ and output space $Y$. 

We show that under stabilizability and detectability assumptions, perturbations of the optimal control problem \eqref{eq:OCPintro} (e.g., of the source term $f\in V^*$) only have an exponentially localized effect. Our results in particular encompass equations where such a localized behavior does not hold for the uncontrolled equation, e.g., Helmholtz equations arising as time-harmonic instance of the wave equation~\cite{EsteMele11,HaucPete22}. Therein, the solutions can be highly oscillatory and, as a consequence, even strongly-localized source terms may have a global effect.

Since such an exponential decay of perturbations yields important insight into localized sensitivity properties of the optimal control problem and reveals robustness of the solution w.r.t.\ perturbations for large spatial and temporal domains, there are various important applications. 
In view of numerical discretization, this locality justifies local refinements of the spatial grid, e.g., on parts of the domain that are of particular interest. 
Other applications are splitting or domain decomposition methods~\cite{LagnLeug04,BartHein06} appearing, e.g., in distributed optimization~\cite{yang2019survey}, where problems on smaller domains are solved and then boundary data on the subdomains is updated iteratively. 
As in the course of the iteration, the boundary conditions at the coupling conditions are inexact, 
analyzing the sensitivity of the solution on the subdomains w.r.t.~these perturbed boundary conditions is crucial~\cite{ShinFaul19,NaShin22}. As further application of the presented analysis, we mention the construction of separable optimal value functions \cite{SperSalu23}. There, a strongly localized sensitivity property of the optimal control problem enables an efficient approximation of the optimal value function by means of neural networks. 

We briefly mention related works considering localized sensitivity in optimal control. In recent years, sensitivity in time was extensively studied. There, the decay of perturbations in the temporal boundary values is manifested in the turnpike property~\cite{FaulGrun22}. 
It relates the dynamic problem to a steady state problem and states that for long time horizons, the solutions of the former are (exponentially) close to the solutions of the latter, see, e.g., \cite{DammGrun14,TrelZuaz15,GugaTrel16,FaulFlas22} and \cite{GuglLi24} for recently provided necessary conditions. The turnpike property was successfully leveraged for time domain decomposition~\cite{ShinFaul19}, multiple shooting~\cite{TrelZuaz15}, neural networks~\cite{EsteGesh20,GeshZuaz22,FaulHemp21}, and efficient mixed-integer programming~\cite{FaulMurr20} and is a central ingredient for performance and stability of model predictive control without terminal ingredients \cite{Grun13,Grun16}. 
Complementary to the decay of the temporal boundary values, another consequence of localized sensitivities in time is the decay of residuals, stemming e.g.~from time or space discretization. 
For linear quadratic optimal control, and under stabilizability and detectability assumptions, such a property was proven in \cite{GrunScha19,GrunScha21} for linear and nonlinear parabolic problems and in~\cite{GrunScha20} in a semigroup context. This decay of discretization errors provides the foundation for efficient a-posteriori goal-oriented discretizations in model predictive control, where only a part on the domain is refined~\cite{GrunScha22} or a-priori coarsening strategies for the time grid~\cite{ShinZava21}. To extend the well-established sensitivity analysis in the (one-dimensional) time to higher dimensions, recent works provide a exponentially localized sensitivity of graph-structured finite-dimensional problems~\cite{ShinAnit22} with promising applications to graph decompositions~\cite{NaShin22,Shin21} and near-optimal and efficient linear quadratic regulators~\cite{ShinLin23}. However, to this date, there is no corresponding result providing an exponential decay of perturbations in space for problems involving PDEs.

This work serves as a starting point in this direction. Under stabilizability and detectability conditions on the pairs $(A,B)$ and $(A,C)$ in~\eqref{eq:OCPintro}, we show that such an exponential decay property holds. 
Our main results are as follows: First, for elliptic problems, and even if the PDEs solution operator does not enjoy exponential locality of perturbations, 
we show that optimal solutions exhibit an exponentially localized sensitivity. We illustrate the developed theory using various examples such as Helmholtz and Poisson equations as well as coupled systems. Second, we provide an extension to parabolic optimal control problems involving advection-diffusion-reaction systems \changed{abstractly given by
	\begin{align*}
	\begin{split}
	\min_{(y,u)\in W(0,T)\times L^2(0,T;U)} \frac12 (\|C(y-y_d)\|_{L^2(0,T;Y)}^2 &+ \|u\|_{L^2(0,T;U)}^2)\\
	\text{s.t. } \dot y + Ay + Bu &= f, \qquad  
	y(0)= y^0.
	\end{split}
	\end{align*}}
\noindent The remainder of this work is organized as follows. In the subsequent Section~\ref{sec:motivation}, we provide a motivation of the presented theory, illustrating the exponential decay in optimal control by means of prototypical one-dimensional elliptic problems. 
In Section~\ref{sec:elliptic}, we prove this property for stationary optimal control problems governed by second-order differential operators and illustrate the results by means of various numerical examples in two dimensions in Section~\ref{sec:num_ell}. Then, in Section~\ref{sec:parab}, we provide an extension to optimal control of parabolic advection-diffusion-reaction systems, which is illustrated in Section~\ref{sec:num_parab} by means of numerical examples. Last, conclusions are drawn in Section~\ref{sec:conclusions}.
\\[.5em]
\noindent \textbf{Notation.} By $\Omega \subset \mathbb{R}^d$, $d\in \mathbb{N}$, we denote a measurable, open and bounded set with Lipschitz boundary. We denote the usual Lebesgue spaces of real-valued and $p$-integrable functions on $\Omega$ by $L^p(\Omega)$, $p\in \mathbb{N}$. 
Correspondingly, we refer to $L^\infty(\Omega)$ as the space of measurable and essentially-bounded functions. 
The Sobolev space of square-integrable and weakly-differentiable functions with square-integrable derivative is denoted by~$H^1(\Omega)$. 
For a measurable subset of the boundary $\Gamma_\mathrm{D}\subset \partial \Omega$, we set $H^1_{\Gamma_\mathrm{D}}(\Omega) = \{v\in H^1(\Omega)\,|\, v=0\ \mathrm{a.e.\ on}\  \Gamma_\mathrm{D}\}$, which is also a Hilbert space when endowed with the $H^1$-topology by continuity of the trace operator $\gamma:H^1(\Omega)\to L^2(\Gamma_\mathrm{D})$. 
If $\Gamma_\mathrm{D}=\partial \Omega$, we write $H^1_0(\Omega) = H^1_{\Gamma_\mathrm{d}}(\Omega)$. To formulate time-dependent problems and for a Banach space $X$, we utilize the Bochner spaces $L^p(0,T;X)$, $p\in \mathbb{N}$, of $p$-integrable $X$-valued functions. 
Correspondingly, $H^1(0,T;X)$ denotes the space of functions in $L^2(0,T;X)$ with weak time derivative in $L^2(0,T;X)$. 
We write $C(\overline{\Omega};X)$, for the space of continuous functions with values in $X$ and abbreviate $C(\overline{\Omega})=C(\overline{\Omega};\mathbb{R})$. 

\section{Motivation: Spatial decay in one dimension}\label{sec:motivation}

\noindent We provide a brief motivation of our analysis by means of three prototypical partial differential equations~(PDEs) on a one-dimensional spatial domain, i.e., the Poisson equation, the screened Poisson equation, and the Helmholtz equation. 
First, we recap well-known properties of these equations to showcase the different stability behavior in the uncontrolled case with a particular focus on the influence of perturbations. 
Then, we illustrate that in optimal control, all equations behave very similarly and enjoy an exponential decay of perturbations. \\[.5em]

\noindent \textbf{Poisson equation with Dirichlet boundary conditions: linear decay.} As a first prototypical example and for domain size~$L$, $L>0$, we consider the one-dimensional elliptic PDE
\begin{align*}
-y''(\omega) = f(\omega)\quad \forall\, \omega \in (0,L) \quad\text{with Dirichlet boundary conditions }\quad y(0)=a,\, y(L)=b
\end{align*}
with $a,b\in \mathbb{R}$ and source term $f:(0,L)\to \mathbb{R}$. 
The solution can be obtained via the linearly decaying Green's function, cf.\ \cite[Section 1.2]{StakHols11},
\begin{align}\label{eq:green}    
\begin{split}
g(\omega,\xi) :=\begin{cases}
(L-\xi)\omega,\qquad 0\leq \omega < \xi,\\
(L-\omega)\xi,\qquad \xi \leq \omega \leq L
\end{cases}
\end{split}
\end{align}
by means of the convolution
\begin{align}\label{eq:conv}
y(\omega) = \int_0^L g(\omega,\xi) f(\xi)\,\mathrm{d}\xi + \changed{a\left(1-\frac{\omega}{L}\right)} + \changed{\frac{b}{L}}\omega.
\end{align}
In this work, we are interested in the stability behavior of the solution~$y$ in $L^p(0,L)$-norms, $1\leq p < \infty$, when perturbing the data $(f,a,b)$, with particular focus on increasing domain sizes $L$. To illustrate this, we first consider a perturbation of the boundary condition on the right. Due to linearity of the problem, we may consider $a=0$ and $f=0$. Then, the sensitivity of the perturbation $b\in \mathbb{R}$ in the $L^2(0,L)$-norm is \changed{increasing} in the domain size, as, in view of \eqref{eq:conv},
$$\|y\|_{L^2(0,L)}^2 = \int_0^L |y(\omega)|^2\,\mathrm{d}\omega = \int_0^L \left|\changed{\frac{b}{L}\omega}\right|^2\,\mathrm{d}\omega = \changed{\frac{b^2}{L^2}}\frac{L^3}{3} = \changed{\frac{b^2L}{3}}.$$
Second, we consider perturbations of the source term~$f$ and thus, again due to linearity, set $a=b=0$. We consider a perturbation $f(\xi) = \delta_{L/2}(\xi)$, with a Dirac delta functional $\delta_{L/2}\in H^{1}(0,L)^*$ modeling a point evaluation at $L/2$, i.e., $\langle \delta_{L/2},g\rangle_{H^1(0,L)^*,H^{1}(0,L)} = h(L/2)$ for $h\in H^1(0,L)\hookrightarrow C([0,L])$. Then $y(\omega) = \int_0^L g(\omega,\xi) \delta_{L/2}(\xi)\,\mathrm{d}\xi = g(\omega,{L/2})$ and  
\begin{align*}
\|y\|_{L^2(0,L)}^2 = \int_0^L |g(\omega,{L/2})|^2\,\mathrm{d}\omega &= \changed{\int_0^{L/2} \left((L-L/2)\omega\right)^2 \,\mathrm{d}\omega + \int_{L/2}^{L} \left((L-\omega)L/2\right)^2 \,\mathrm{d}\omega} \\
&= \changed{\frac{L^2}{4} \left(\int_0^{L/2} \omega^2\,\mathrm{d}\omega + \int_{L/2}^L (L-\omega)^2\,\mathrm{d}\omega\right) = \frac{L^2}{4} \cdot \frac{L^3}{12}  = \frac{L^5}{48}.}
\end{align*}
Thus, we may conclude for this Poisson equation that the influence of the perturbations, and thus also the norm of the linear solution operator measured in integral norms, is not uniformly bounded in the domain size $L$. 
\\[.5em]
\noindent \textbf{Screened Poisson equation: exponential decay.} As a second example, we let $k \in \mathbb{R}\setminus \{0\}$ and consider the screened Poisson equation
\begin{align*}
-y''(\omega) + k^2 y(\omega) &= f(\omega) \qquad \forall \omega \in (0,L)
\end{align*}
with $f:(0,L)\to \mathbb{R}$ and either Dirichlet or Neumann boundary conditions, i.e.,
\begin{align*}
y(0)=a,\, y(L)=b \qquad \mathrm{or} \qquad     y'(0)=a,\, y'(L)=b
\end{align*}
for $a,b\in \mathbb{R}$.
Here, we directly see that the solution operator norm is bounded uniformly in the domain size. 
For example, 
for $a=b=0$ and a test function $v\in H_0^1(0,L)$ (in case of Dirichlet conditions) or $v\in H^1(0,L)$ (in case of Neumann boundary conditions), integration by parts yields
\begin{align*}
\langle y',v'\rangle_{L^2(0,L)} + k^2\langle y,v\rangle_{L^2(0,L)} = \langle f,v\rangle_{H^1(0,L)^*,H^1(0,L)}.
\end{align*}
Thus, choosing the solution $y$ as a test function, i.e., $v=y$, we get
\begin{align*}
\min\{1,k^2\} \|y\|_{H^1(0,L)} \leq \|f\|_{H^1(0,L)^*}.
\end{align*}
\changed{As $k$ is independent of the domain size and if the perturbation $f$ is bounded uniformly in $L$, this yields an upper bound uniform in $L$. This, in turn implies that $\|y\|^2_{H^1(0,L)}$ is bounded uniformly in the domain size, that is, for increasing domain sizes $L$, the state~$y$ and its weak gradient~$\nabla y$ have to be small on the majority of the domain. More precisely and e.g.\ for the state $y$, for any $\varepsilon>0$ and by Markov's inequality we have
	\begin{align*}
	\operatorname{meas}\left(\{\omega \in [0,L]\,|\, |y(t)|>\varepsilon\}  \right) \leq \frac{1}{\varepsilon^2} \|y\|_{L^2(0,L)}^2 \leq \frac{\|f\|^2_{H^1(0,L)^*}}{\min\{1,k^4\}\varepsilon^2}
	\end{align*}
	with an upper bound that is uniform in $L$.}\\[.5em]
\noindent \textbf{Helmholtz Problem: No decay.} Last, for a wave number $k\in \mathbb{R} \setminus \{ 0 \}$, we consider the Helmholtz problem 
\begin{align*}
-y''(\omega) - k^2y(\omega) = f(\omega)\quad \forall\, \omega \in (0,L) \quad\text{with Dirichlet boundary conditions}\quad y(0)=a, y(L)=b
\end{align*}
for $a,b\in \mathbb{R}$ and $f:(0,L) \to \mathbb{R}$. This Helmholtz equation is known to exhibit oscillatory solution behavior as it occurs, e.g., as the time-harmonic case of the wave equation~\cite{EsteMele11,HaucPete22}. Further, it can also occur in semi-discretizations of a reaction-diffusion equation, i.e.,
\begin{align*}
\partial_ty(t,\omega) = \partial_{\omega\omega} y(t,\omega) + c^2 y(t,\omega)
\end{align*}
with $c\in \mathbb{R}$, as, using the implicit Euler method with stepsize $\Delta t >0$ yields
\begin{align*}
(I + \Delta t\cdot(-\partial_{\omega\omega} -c^2))y(t_{k+1},\cdot) = y(t_{k},\cdot).
\end{align*}
If we divide this equation by $\Delta t$, the left-hand side is given by $-\partial_{\omega\omega}y + (\frac{1}{\Delta t} - c^2)y$.
If the time step $\Delta t$ is not small enough w.r.t.~the level of instability $c^2$, this leads to Helmholtz-type problems with \changed{$k^2 =(c^2- \frac{1}{\Delta t})$}.\\[.5em]

\begin{figure}[H]
	\centering
	\includegraphics[width=.45\linewidth]{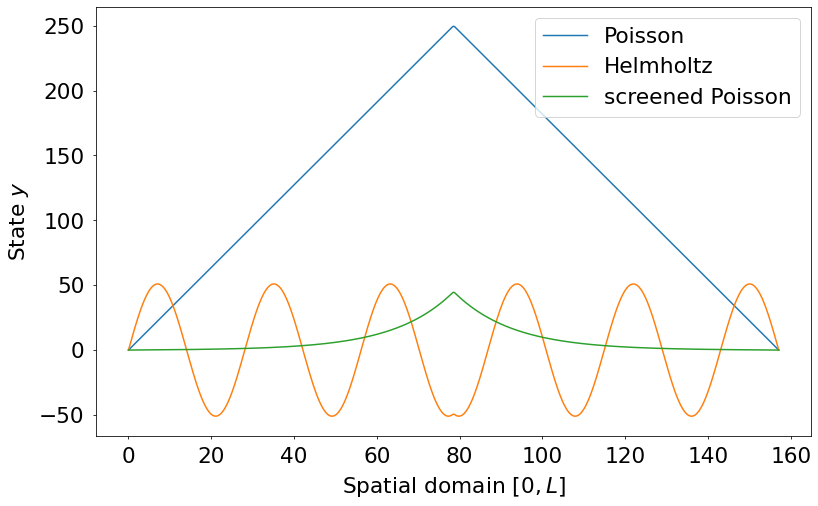}
	\includegraphics[width=.462\linewidth]{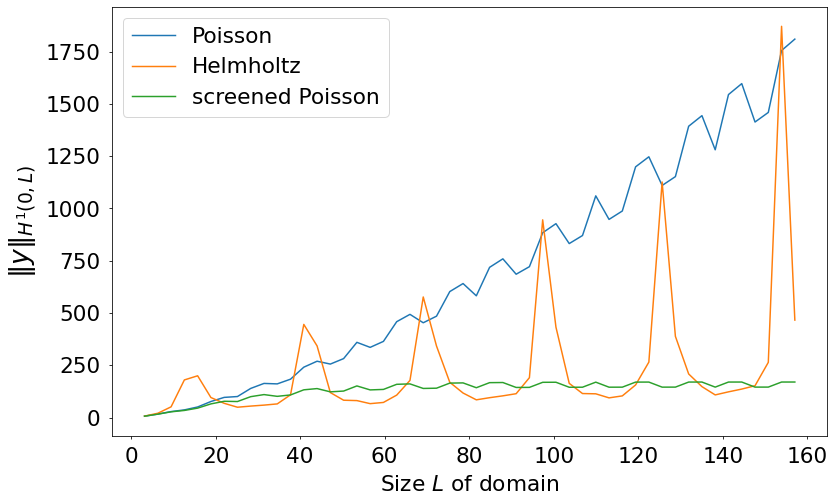}\\
	\includegraphics[width=.45\linewidth]{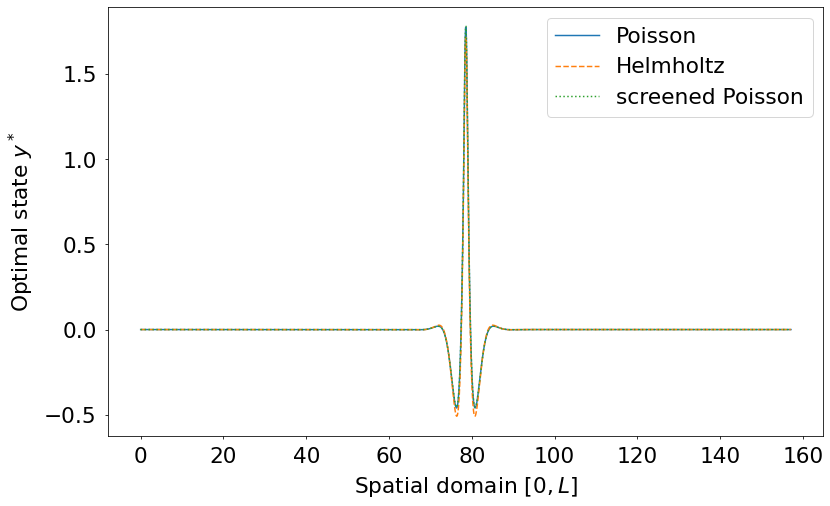}\hspace{.35em}
	\includegraphics[width=.45\linewidth]{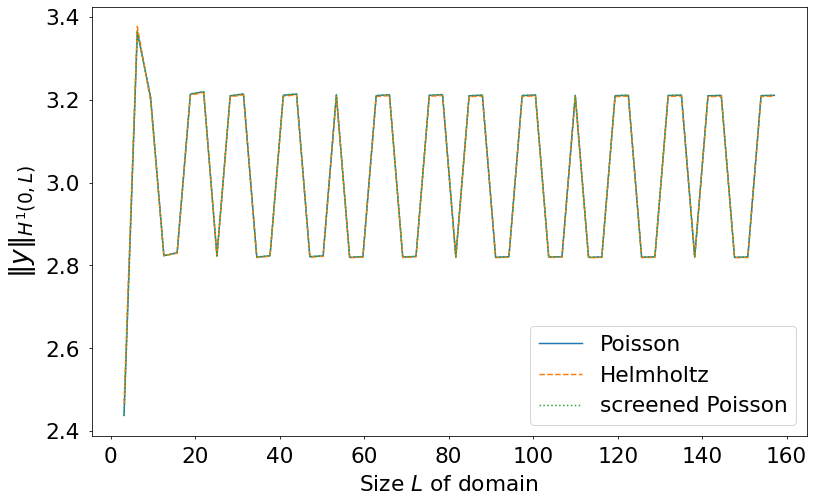}
	\caption{Behavior of the three PDEs for a source term localized in the middle of the domain. Top: Simulation, Bottom: Optimal control.}
	\label{fig:motivation}
\end{figure}

\noindent \textbf{Optimal control: Exponential decay.} As motivated by theoretical considerations above, these three classes of stationary equations have a very different sensitivity w.r.t.\ perturbations. This behavior is numerically illustrated in the upper row of Figure~\ref{fig:motivation} for a localized, piecewise constant perturbation $f$ centered in the middle of the domain. 
We see in the top left plot that only the screened Poisson equation exhibits a localized effect of the perturbation, whereas the influence for the Poisson equation is linearly localized (due to linear decay of the Green's function \eqref{eq:green} in one spatial dimension) and for the Helmholtz equation it is global (due to a lack of coercivity). 
Correspondingly, the norm of the solution increases for larger domain sizes, as depicted in the top right plot of Figure~\ref{fig:motivation}. 
In the bottom row, we depict 
the optimal solutions with distributed control and observation, i.e., setting $C=B=I_{L^2(0,L)}$ in~\eqref{eq:OCPintro}. 
We observe that despite the very different behavior in the uncontrolled case, the optimal states show a strongly (exponentially) localized behavior. As a consequence, the solution norms are all bounded uniformly in the domain size, as depicted in the bottom right plot of Figure~\ref{fig:motivation}. \changed{For details regarding the implementation, we refer to Section~\ref{sec:num_ell}.}

\noindent In this work, we will rigorously show that -- regardless of the stability of the uncontrolled equation -- this exponentially localized effect of perturbations of the optimal control problem~\eqref{eq:OCP} can be guaranteed under stabilizability and detectability assumptions on the involved operators.

\section{Spatial decay of perturbations: Elliptic optimal control}\label{sec:elliptic}

\noindent In the following, we always assume a measurable, open and bounded domain $\Omega\subset \mathbb {R}^d$, $d\in \mathbb{N}$, with Lipschitz boundary $\partial \Omega$. 
Further, let $U$ and $Y$ be Hilbert spaces, which will be tacitly identified with their respective dual spaces. 
Moreover, $(V,\|\cdot\|_V)$ is assumed to be a reflexive Banach space and will be chosen as a subspace of $H^1(\Omega)$ including suitable boundary conditions depending on the particular application.\\ 
In this section, we present the main result of this work for elliptic optimal control problems of the form
\begin{align}\label{eq:OCP}
\min_{(y,u)\in V\times U} \frac12 (\|C(y-y_d)\|_Y^2 &+ \|u\|_U^2) \qquad\text{subject to}\qquad Ay + Bu = f,
\end{align}
where $B\in L(U,V^*)$, $C\in L(V,Y)$, $f\in V^*$, and $y_d\in V$. 
Here, $A:V\to V^*$ is a surjective (not necessarily boundedly invertible) differential operator on $\Omega$. 
Whereas the part of this subsection is independent of the particular structure of $A$, we will later assume $A$ to be of the form 
$$
\langle Ay,v\rangle_{V^*,V} = \langle a\nabla y,\nabla v\rangle_{L^2(\Omega)} + \langle b^\top\nabla y,v \rangle_{L^2(\Omega)} + \langle cy,v\rangle_{L^2(\Omega)}
$$ with coefficients $a\in L^\infty(\Omega;\mathbb{R}^{d\times d})$ and $b\in L^\infty(\Omega,\mathbb{R}^d)$ and $c\in L^\infty(\Omega)$.\\[.5em]
\noindent A central tool in our analysis will be the first-order optimality system, which --~due to the linear-quadratic nature of the problem~-- poses the necessary and sufficient conditions of optimality. As $A:V\to V^*$ is surjective, $[A,B]:V\times U\to V^*$ is as well and has, hence, in particular closed range. 
Thus, we may use the following standard result on optimality conditions, cf.~\cite[Theorem 1.1 and Remark 1.2]{Schi13}.
\begin{lemma}
	Let $(y^*,u^*)\in V\times U$ be an optimal solution of \eqref{eq:OCP}. Then there is an adjoint state $p\in V$ such that
	\begin{align}\label{eq:optsys}
	\underbrace{\begin{pmatrix}
		C^*C & 0 & A^* \\
		0&I&B^*\\
		A&B&0
		\end{pmatrix}}_{=:M}
	\begin{pmatrix}
	y^*\\u^*\\p
	\end{pmatrix}
	= \begin{pmatrix}
	C^*Cy_d\\0\\f
	\end{pmatrix}.
	\end{align}
\end{lemma}

\noindent Based on this optimality condition characterizing optimal solutions, we will prove a result showing locality of the influence of perturbations of the right-hand side. 
Therein, it will be crucial to have a bound of the solution operator, which is uniform in the size of the spatial domain $\Omega$. 
Hence, before providing the main result of this work in Theorem~\ref{thm:scaling} of Subsection~\ref{subsec:scaling}, we prove such a uniform bound under stabilizability assumptions in Theorem~\ref{thm:opnorm} of the following subsection. The approach is closely related to the sensitivity analysis for temporal decay of perturbations presented in~\cite{GrunScha19}, where uniform bounds in terms of the time horizon were central.

\subsection{Bounds on the solution operator --~uniform in the domain size}\label{subsec:solnorm}

\noindent The following theorem shows that, under uniform stabilizability and detectability assumptions, the solution operator of the optimality system, and hence of the optimal control problem, is uniformly bounded in the domain size. 
It is motivated by the prototypical examples depicted in Figure~\ref{fig:motivation}, where a uniform bound on the solution was an indicator for an exponential decay of the perturbation. 
This uniform bound will play a central role in the main result of this section stated in Theorem~\ref{thm:scaling}, in which we prove an exponential decay of perturbations.
\begin{theorem}\label{thm:opnorm}
	Let the operators $A \in L(V,V^*)$, $B \in L(U,V^*)$, and $C \in L(Y,V^*)$ be bounded uniformly in~$| \Omega |$. 
	Suppose the existence of a constant~$\alpha > 0$ 
	independent of the domain size~$|\Omega|$ and operators $K_B\in L(V,U)$, $K_C\in L(Y,V^*)$, which are bounded uniformly in~$| \Omega |$, such that the stabilizability and detectability conditions
	\begin{enumerate}
		\item[(i)] 
		$\langle (A+BK_B)y,y\rangle_{V^*, V} \geq \alpha \|y\|_{V}^2$ and
		\item[(ii)] 
		$\langle (A+K_CC)y,y\rangle_{V^*, V} \geq \alpha \|y\|_{V}^2$
	\end{enumerate}
	hold.
	Then there exists a constant $c \geq 0$ independent of $|\Omega|$ such that
	\begin{align*}
	\changed{\|M^{-1}\|_{L(V^*\times U\times V^*,V\times U\times V )}} &\leq c.
	\end{align*}
	for the operator matrix~$M$ defined in the optimality system~\eqref{eq:optsys}.
\end{theorem}
\begin{proof}
	\changed{Let $(y,u,p) \in V\times U \times V$ solve \eqref{eq:optsys} with right-hand side $(\varepsilon_1,\varepsilon_2,\varepsilon_3)\in V^* \times U \times V^*$. We will now prove $\|(y,u,p)\|_{V\times U\times V} \leq c \|(\varepsilon_1,\varepsilon_2,\varepsilon_3)\|_{V^*\times U\times V^*}$ with $c\geq 0$ independent of $|\Omega|$ which implies the claimed bound on the solution operator $M^{-1}$.
		For brevity of notation, we will denote $ (y,v) := \langle y,v\rangle_{V^*, V}$ and $\langle y,v\rangle = \langle y,v\rangle_{L^2(\Omega)}$. First, we eliminate the control by means of the second equation, i.e.,
		\begin{align}\label{eq:subsU}
		u = -B^* p + \varepsilon_2.
		\end{align}
		Inserting this into \eqref{eq:optsys} with right-hand side $(\varepsilon_1,\varepsilon_2,\varepsilon_3)$, we get
		\begin{align}\label{eq:optcond_elim}
		\begin{pmatrix}
		C^*C & A^*\\
		A & -BB^*
		\end{pmatrix}
		\begin{pmatrix}
		y\\p
		\end{pmatrix}
		=
		\begin{pmatrix}
		\varepsilon_1\\
		\varepsilon_3 - B \varepsilon_2
		\end{pmatrix}.
		\end{align}
		Let $K_C \in L(Y,V^*)$ as in Assumption~(ii). Then, reformulating the state equation and adding $(K_CCy,v)$ on both sides, we get
		\begin{align*}
		((A+K_CC)y,v) = (\varepsilon_3-B\varepsilon_2 + BB^*p + K_CCy,v)
		\end{align*}
		for all $v\in V$. Thus, testing this equation with $y$ and invoking Assumption~(i), we get
		\begin{align*}
		\alpha \|y\|_V^2 &\leq (\varepsilon_3-B\varepsilon_2 + BB^*p + K_CCy,y) 
		\\&\leq \|\varepsilon_3-B\varepsilon_2\|_{V^*}\|y\|_{V} + \|B\|_{L(U,V^*)}\|B^*p\|_{U}\|y\|_V + \|K_C\|_{L(Y,V^*)}\|Cy\|_Y \|y\|_V.
		\end{align*}
		We divide this estimate by $\alpha \|y\|_V$ and deduce
		\begin{align*}
		\|y\|_V \leq \frac{1}{\alpha}\left( \|\varepsilon_3\|_{V^*} +\|B\|_{L(U,V^*)}\|\varepsilon_2\|_U+ \max\{\|B\|_{L(U,V)}, \|K_C\|_{L(Y,V)}\}\left(\|B^*p\|_{U}+\|Cy\|_Y\right)\right).
		\end{align*}
		Squaring this bound and applying $(a+b)^2\leq 2(a^2+b^2)$ for $a,b\in \mathbb{R}$ three times yields
		\begin{align}\label{eq:statebound}
		\|y\|^2_V \leq \changed{\frac{4}{\alpha^2}\left( \|\varepsilon_3\|^2_{V^*} +\|B\|^2_{L(U,V^*)}\|\varepsilon_2\|^2_U+ \max\{\|B\|_{L(U,V)}, \|K_C\|_{L(Y,V)}\}^2\left(\|B^*p\|^2_{U}+\|Cy\|^2_Y\right)\right)}.
		\end{align}
		We now prove a similar bound for the adjoint equation
		\begin{align*}
		(C^*Cy + A^*p,v) = (\varepsilon_1, v)
		\end{align*}
		for all $v\in V$. To this end, let $K_B\in L(V,U)$ as in Assumption~(i) such that
		\begin{align*}
		((A+BK_B)p,v) = (\varepsilon_1 - C^*Cy + BK_Bp,v)
		\end{align*}
		for all $v\in V$. Testing this equation with the adjoint state $p$ and proceeding analogously as before for the state yields
		\begin{align}\label{eq:adjointbound}
		\|p\|^2_V \leq \frac{2}{\alpha^2}\left(\|\varepsilon_1\|^2_{V^*} + 2\max\{\|C^*\|_{L(Y,V^*)},\|K_B\|_{L(V,U)}\}^2 \left(\|Cy\|^2_Y + \|B^*p\|^2_U\right)\right).
		\end{align}
		We now remain to bound the term $\|Cy\|^2_Y + \|B^*p\|^2_U$ by means of the norms of $\varepsilon_i$, $i=1,2,3$. To this end, we test the adjoint equation, i.e., the first equation of \eqref{eq:optcond_elim} with the state, and the state equation, i.e., the second equation of \eqref{eq:optcond_elim} with the adjoint and subtract the results, leading to
		\begin{align*}
		(C^*Cy,y) + (A^*p,y) - (Ay,p) + (BB^*p,p) = (\varepsilon_1,y) - (\varepsilon_3,p) + (B\varepsilon_2,p),
		\end{align*}
		which implies
		\begin{align}\label{eq:costbound}
		\|Cy\|^2_Y + \|B^*p\|_U^2 \leq \|\varepsilon_1\|_{V^*}\|y\|_V + \left(\|\varepsilon_3\|_{V^*} + \|B\|_{L(U,V^*)}\|\varepsilon_2\|_U\right)\|p\|_V.
		\end{align}
		Adding the two inequalities \eqref{eq:statebound} and \eqref{eq:adjointbound} and plugging in \eqref{eq:costbound}, we obtain
		\begin{align}\label{eq:sumbound}
		\begin{split}
		\|y\|^2_V+\|p\|^2_V \leq 
		&\changed{\frac{2}{\alpha ^2} \left( \|\varepsilon _1\|_{V^*} ^2 + 2 \|B\|_{L(U,V^*)}^2 \|\varepsilon_2\|_{U}^2+2\|\varepsilon_3\|_{V^*} ^2 \right)} \\&+ \frac{2d}{\alpha^2}\left( \|\varepsilon_1\|_{V^*}\|y\|_V + \left(\|\varepsilon_3\|_{V^*} + \|B\|_{L(U,V^*)}\|\varepsilon_2\|_U\right)\|p\|_V\right)
		\end{split}
		\end{align}
		with $d =2\left(\max\{\|B\|_{L(U,V)}, \|K_C\|_{L(Y,V)}\}^2+\max\{\|C^*\|_{L(Y,V^*)},\|K_B\|_{L(V,U)}\}^2\right)$. We now apply the estimate $ab\leq \tfrac12(r^2 a^2 + \frac{1}{r^2}b^2)$ for all $a,b,r\in \mathbb{R}$ twice, i.e., with
		\begin{align*}
		a = \|y\|_V,\ b = \|\varepsilon_1\|_{V^*},\ r = \sqrt{\frac{\alpha^2}{2d}}  \quad \mathrm{and}\quad a = \|p\|_V,\ b = \|\varepsilon_3\|_{V^*} + \|B\|_{L(U,V^*)}\|\varepsilon_2\|_U,\ r = \sqrt{\frac{\alpha^2}{2d}},
		\end{align*}
		respectively, to deduce
		\begin{align*}
		\frac{2d}{\alpha^2}&\left( \|\varepsilon_1\|_{V^*}\|y\|_V + \left(\|\varepsilon_3\|_{V^*} + \|B\|_{L(U,V^*)}\|\varepsilon_2\|_U\right)\|p\|_V\right) 
		\\&\leq \frac{2d^2}{\alpha^4}\left( \|\varepsilon_1\|^2_{V^*} \right.\left.+ \left(\|\varepsilon_3\|_{V^*} + \|B\|_{L(U,V^*)}\|\varepsilon_2\|_U\right)^2\right) + \frac12\left(\|y\|_V^2 + \|p\|_V^2 \right).
		\end{align*}
		Inserting this estimate into \eqref{eq:sumbound} and subtracting $\frac12\left(\|y\|_V^2 + \|p\|_V^2 \right)$ yields
		\begin{align*}
		\frac12\left(\|y\|_V^2 + \|p\|_V^2 \right) &\leq \frac{2}{\alpha^2}\left( \|\varepsilon_1\|^2_{V^*} + 2\|B\|_{L(U,V^*)}\|\varepsilon_2\|_U + \|\varepsilon_3\|_{V^*}^2\right)\\
		& \qquad + \frac{2d^2}{\alpha^4}\left(\|\varepsilon_1\|^2_{V^*} + (\|\varepsilon_3\|_{V^*} + \|B\|_{L(U,V^*)}\|\varepsilon_2\|_U)^2 \right)\\
		&\leq \frac{2}{\alpha^2}\left( \|\varepsilon_1\|^2_{V^*} + 2\|B\|_{L(U,V^*)}\|\varepsilon_2\|_U + \|\varepsilon_3\|_{V^*}^2\right)\\
		& \qquad + \frac{2d^2}{\alpha^4}\left(\|\varepsilon_1\|^2_{V^*} + 2\|\varepsilon_3\|_{V^*}^2 + 2\|B\|_{L(U,V^*})^2\|\varepsilon_2\|_U^2 \right)
		\\&\leq \frac{2}{\alpha^2}\left(1+\frac{d^2}{\alpha^2}\right)\left( \|\varepsilon_1\|^2_{V^*} + 2\|\varepsilon_3\|_{V^*}^2 + 2\|B\|^2_{L(U,V^*)}\|\varepsilon_2\|^2_U\right)\\
		&\leq 4\left(\frac{1}{\alpha^2}+d^2 \right)\max\{1,\|B\|^2_{L(U,V^*)}\}\left(\|\varepsilon_1\|^2_{V^*} + \|\varepsilon_2\|^2_{U} + \|\varepsilon_3\|^2_{V^*}\right).
		\end{align*}
		Hence,
		\begin{align*}
		\|y\|^2_V + \|p\|^2_V \leq \tilde c \left( \|\varepsilon_1\|^2_{V^*} +  \|\varepsilon_2\|^2_U +  \|\varepsilon_3\|^2_{V^*} \right)
		\end{align*}
		with $\tilde c :=8 \left(\frac{1}{\alpha^2} + d^2\right)\max\{1,\|B\|_{L(U,V^*)}\}$. 
		Last, the bound for the control follows by resubstituting \eqref{eq:subsU}:
		\begin{align*}
		\|u\|^2_U = \|-B^* p + \varepsilon_2\|^2_U \leq 2\left(\|B\|_{L(U,V^*)}^2 \|p\|^2_V + \|\varepsilon_2\|^2_U\right)
		\end{align*}
		such that the claim follows with $c = \tilde c + 2(\|B\|_{L(U,V^*)}\tilde c +1)$.}
\end{proof}

\noindent We note that the proof of Theorem~\ref{thm:opnorm} \changed{shows that $M$ is injective\footnote{\changed{Surjectivity may also be straightforwardly shown due to existence of optimal solutions for any datum, and the fact that any optimal solution solves the optimality system. More precisely, for a given right-hand side $(\ell_1,\ell_2,\ell_3)$ of \eqref{eq:optsys}, we may define an optimal control problem with cost functional $\|Cy\|_Y^2 + \|u\|_U^2 + \langle \ell_1,y\rangle_V + \langle \ell_2,u\rangle_U$ and constraint $Ay+Bu=\ell_3$. The corresponding optimal triple $(y^*,u^*,p)$, which exists due to standard argumentation, solves \eqref{eq:optsys} with $(\ell_1,\ell_2,\ell_3)$ as right-hand side, i.e., $M$ is surjective.}}, that is, if $(\varepsilon_1,\varepsilon_2,\varepsilon_3)=0$, then $(y^*,u^*,p)=0$. Thus,} there is always a unique triple $(y^*,u^*,p)\in V\times U\times V$ solving \eqref{eq:optsys}. This in particular yields uniqueness of the optimal state, despite the fact that the PDE can posses multiple solutions as $A$ is not assumed to be injective.




\subsection{Exponential decay of perturbations}\label{subsec:scaling}

\noindent First, we prove a result for scalings including the one-norm, which we denote by $\|z\|_1 = \sum_{i=1}^d |z_i|$ for $z\in \mathbb{R}^d$. \changed{Further, we denote the vector-valued componentwise sign function by $\operatorname{sgn}(\omega) = (\operatorname{sgn}(\omega_1),\ldots,\operatorname{sgn}(\omega_d))^\top$.} \changed{The result is well-known and straightforwardly follows from the chain rule for weak derivatives \cite[Proposition 9.5]{Brezi11}. We present it here for completeness.}

\begin{lemma}\label{lem:onenorm}
	If $v \in H^1(\Omega)$, then $|v| \in H^1(\Omega)$ with $\nabla |v| = \operatorname{sgn}(v)\nabla v$. Further, if $\rho(\omega) = e^{\mu \|z-\omega\|_1}$ for $z\in \Omega$, then 
	$$\nabla \rho(\omega) = -\mu \begin{pmatrix}
	\operatorname{sgn}(z_1-\omega_1)\\
	\vdots\\
	\operatorname{sgn}(z_d-\omega_d)
	\end{pmatrix} \rho(\omega).$$
\end{lemma}
\begin{proof}
	Clearly, the absolute value function in one space dimension $\varphi(\omega) = |\omega|$ is weakly differentiable with $ \varphi'(\omega) = \operatorname{sgn}(\omega)$. Moreover, $\varphi$ is essentially bounded and the derivative is square essentially bounded as well such that $\varphi \in W^{1,\infty}(\Omega)$. Thus, the first claim follows by the chain rule for weak derivatives. The second claim follows by straightforward computations.
\end{proof}
\noindent Whereas the abstract result of Theorem~\ref{thm:opnorm} did not hinge on a particular structure of $A$, we will consider PDEs governed by second-order differential operators in the following.

\begin{assumption}\label{as:secondorder}
	Let $V\subset H^1(\Omega)$ endowed with the $H^1(\Omega)$-topology. The operator $A:V \to V^*$ is assumed to be a second order differential operator of the form $$\langle Ay,v\rangle_{V^*,V} = \langle a\nabla y,\nabla v\rangle_{L^2(\Omega)} + \langle b^\top\nabla y,v \rangle_{L^2(\Omega)} + \langle cy,v\rangle_{L^2(\Omega)}$$ with coefficients $a\in L^\infty(\Omega;\mathbb{R}^{d\times d})$, $b\in L^\infty(\Omega,\mathbb{R}^d)$ and $c\in L^\infty(\Omega)$. 
\end{assumption}
A central tool to derive a localized sensitivity analysis will be an exponentially scaled version of the optimality system \eqref{eq:optsys}. To this end, we utilize the following auxiliary result.
\begin{lemma}\label{lem:perturbations}
	Let Assumption~\ref{as:secondorder} hold and $y\in V$ solve
	\begin{align*}
	\langle Ay,v\rangle_{V^*,V} = \langle f,v\rangle_{V^*,V} \qquad \forall v\in V
	\end{align*}
	with $f\in V^*$ satisfying $\langle f,\rho y\rangle_{V^*,V} = \langle \rho f,y\rangle_{V^*,V}$ for all $\rho \in W^{1,\infty}(\Omega)$.
	Then the scaled variable $\tilde y(\omega) = e^{\mu\|z-\omega\|_1}y(\omega)$, $z\in \Omega$, $\mu \in \mathbb{R}$, solves
	\begin{align*}
	\langle (A + \mu F_1 - \mu^2 F_2)\tilde y,v\rangle_{V^*,V} = \langle \tilde f,v\rangle_{V^*,V} \qquad \forall v\in V
	\end{align*}
	with scaled right-hand side $\tilde f(\omega) = e^{\mu \|z-\omega\|_1}f(\omega)$ and 
	operators $F_1,F_2:V\to V^*$ defined by
	\begin{align}\label{eq:defnF_1norm}
	\begin{split}
	\langle F_1 y, v\rangle_{V^*,V} &=\langle a\operatorname{sgn}(z-\cdot) y, \nabla v\rangle - \langle a\nabla y, \operatorname{sgn}(z-\cdot) v\rangle + \langle b^\top\operatorname{sgn}(z-\cdot) y,v\rangle\\
	\mathrm{and} \qquad  \langle F_2 y, v\rangle_{V^*,V} &= \langle a\operatorname{sgn}(z-\cdot)y,\operatorname{sgn}(z-\cdot)v\rangle.  
	\end{split}
	\end{align}
	Moreover,
	\begin{align}\label{eq:F_bdd_1norm}
	\|F_1\|_{L(V,V^*)}\leq \max\{\|a\|_{L^\infty(\Omega)},\|b\|_{L^\infty(\Omega;\mathbb{R}^d)}\} \qquad \mathrm{and}\qquad \|F_2\|_{L(V,V^*)}\leq \|a\|_{L^\infty(\Omega)}.
	\end{align}
\end{lemma}
\begin{proof}
	We set $\rho(\omega) = e^{\mu \|z-\omega\|_1}$. By means of the product rule and Lemma~\ref{lem:onenorm}, we compute
	\begin{align*}
	\nabla \tilde y(\omega)=  \nabla \rho(\omega) y(\omega) + \rho(\omega) \nabla y(\omega) = -\mu \operatorname{sgn}(z-\omega) \rho(\omega) y(\omega) + \rho(\omega)\nabla y(\omega).
	\end{align*}
	Thus,
	\begin{align*}
	\langle a\nabla y,\nabla v\rangle &= \langle a\nabla \tilde y + a\mu \operatorname{sgn}(z-\cdot)\tilde y,\rho^{-1}\nabla v \rangle.
	\end{align*}
	We want to pass to test functions $\tilde v = \rho^{-1} v =  e^{-\mu\|z-\cdot\|_1}v$ and compute
	\begin{align*}
	\rho^{-1} \nabla v =  \nabla (\rho^{-1} v) - \mu \operatorname{sgn}(z-\cdot) \rho^{-1} v
	\end{align*}
	such that
	\begin{align*}
	\langle a\nabla \tilde y + a\mu \operatorname{sgn}(z-\cdot) \tilde y,\rho^{-1}\nabla v \rangle = \langle a\nabla\tilde y +a\mu \operatorname{sgn}(z-\cdot) \tilde y, \nabla \tilde v - \mu \operatorname{sgn}(z-\cdot)\tilde v \rangle .
	\end{align*}
	Further, as
	\begin{align*}
	\quad \langle b^\top \nabla y,v\rangle &= \langle b^\top\nabla \tilde y + b^\top \mu \operatorname{sgn}(z-\cdot) \tilde y,\rho^{-1} v\rangle = \langle b^\top\nabla \tilde y + b^\top \mu \operatorname{sgn}(z-\cdot) \tilde y,\tilde v\rangle\\
	\quad \mathrm{and} \quad \langle cy,v\rangle &= \langle c\tilde y, \rho^{-1} v\rangle = \langle c\tilde y,\tilde v\rangle
	\end{align*}
	we may rewrite
	\begin{align*}
	\langle Ay,v\rangle_{V^*,V} &=\langle a\nabla y,\nabla v\rangle + \langle b^\top \nabla y,v\rangle+ \langle cy,v\rangle \\
	&= \langle a\nabla \tilde y\!+\! a\mu \operatorname{sgn}(z-\cdot) \tilde y, \nabla \tilde v \!-\! \mu \operatorname{sgn}(z-\cdot)\tilde v \rangle \!+\! \langle b^\top\nabla \tilde y \!+\! b^\top \mu \operatorname{sgn}(z-\cdot) \tilde y,\tilde v\rangle\!+\! \langle c\tilde y,\tilde v\rangle\\
	&= \langle A \tilde y,\tilde v\rangle_{V^*,V} + \mu\left(\langle a\operatorname{sgn}(z-\cdot) \tilde y, \nabla \tilde v\rangle - \langle a\nabla \tilde y, \operatorname{sgn}(z-\cdot) \tilde v\rangle + \langle b^\top\operatorname{sgn}(z-\cdot) \tilde y,\tilde v\rangle\right) \\&
	\qquad \qquad - \mu^2\langle a\operatorname{sgn}(z-\cdot)\tilde y,\operatorname{sgn}(z-\cdot)\tilde v\rangle\\
	&= \langle (A + \mu F_1 - \mu^2 F_2)\tilde y,\tilde v\rangle_{V^*,V}.
	\end{align*}
	Further, for the right-hand side, we get
	\begin{align*}
	\langle f,v \rangle_{V^*,V} =  \langle \rho f, \rho^{-1} v \rangle_{V^*,V} =  \langle \tilde f,\tilde v \rangle_{V^*,V}.
	\end{align*}
	Now the result follows as testing with all $v\in V$ equivalent to testing with all $\tilde v = \rho^{-1} v \in V$.
	
	To prove the bounds \eqref{eq:F_bdd_1norm}, we have $|\operatorname{sgn}(z-\cdot)| =1$ and compute
	\begin{align*}
	&\langle F_1 y, v\rangle_{V^*,V} \\
	&=\langle a\operatorname{sgn}(z-\cdot) y, \nabla v\rangle - \langle a\nabla y, \operatorname{sgn}(z-\cdot) v\rangle + \langle b^\top\operatorname{sgn}(z-\cdot) y,v\rangle \\
	&\leq \max\{\|a\|_{L^\infty(\Omega)},\|b\|_{L^\infty(\Omega;\mathbb{R}^d)}\}\left(\|y\|_{L^2(\Omega)}\|\nabla v\|_{L^2(\Omega;\mathbb{R}^d)} \right.\\&\hspace{6cm}+\left. \|\nabla y\|_{L^2(\Omega;\mathbb{R}^d)}\|v\|_{L^2(\Omega)} \!+\! \|y\|_{L^2(\Omega)}\|v\|_{L^2(\Omega)}\right)\\
	&\leq \frac{\max\{\|a\|_{L^\infty(\Omega)},\|b\|_{L^\infty(\Omega;\mathbb{R}^d)}\}}{2}\left(\|\nabla y\|_{L^2(\Omega;\mathbb{R}^d)}^2 \!+\!\|v\|^2_{L^2(\Omega)} \!+\! \|y\|_{L^2(\Omega)}^2 \right.\\&\left.\hspace{6cm}+ \|\nabla v\|^2_{L^2(\Omega;\mathbb{R}^d)} \!+\! \|y\|^2_{L^2(\Omega)} \!+\! \|v\|^2_{L^2(\Omega)}\right)\\
	&\leq \max\{\|a\|_{L^\infty(\Omega)},\|b\|_{L^\infty(\Omega;\mathbb{R}^d)}\}(\|y\|_{H^1(\Omega)}^2 + \|v\|_{H^1(\Omega)}^2).
	\end{align*}
	Thus, as $\|\cdot\|_V = \|\cdot\|_{H^1(\Omega)}$, we have
	 $$\|F_1\|_{L(V^*,V)} = \sup_{\|y\|_{V}=1}\sup_{\|v\|_{V}=1}\langle F_1 y, v\rangle_{V^*,V} \leq \max\{\|a\|_{L^\infty(\Omega)},\|b\|_{L^\infty(\Omega;\mathbb{R}^d)}\}$$ and the bound for $F_2$ follows directly from
	\begin{align*}
	\langle F_2 y, v\rangle_{V^*,V} = \langle a\operatorname{sgn}(z-\cdot)y,v\operatorname{sgn}(z-\cdot)\rangle &\leq \|a\|_{L^\infty(\Omega)}\|y\|_{L^2(\Omega)}\|v\|_{L^2(\Omega)}
	\\&\leq \|a\|_{L^\infty(\Omega)}\|y\|_{H^1(\Omega)}\|v\|_{H^1(\Omega)}.
	\end{align*}
\end{proof}
\noindent  We note that the assumption on the right-hand side in above Lemma~\ref{lem:perturbations} is crucial to shift the scaling from the test functions to the functional and may be violated for non-locally acting linear functionals, e.g., given by convolutions. However, it is satisfied for a large class of right-hand sides appearing in applications, e.g., if $\langle f,y\rangle_{V^*,V} = \int_\Omega f(\omega)y(\omega)\,\mathrm{d}\omega$. A similar assumption is also necessary for the input and the output operator as stated in the following.

\begin{assumption}\label{as:local}
	The input and output operators \changed{satisfy the following: For any} $\rho\in W^{1,\infty}(\Omega)$, we have
	\begin{align*}
	\langle B(\rho u),v\rangle_{V^*,V} = \langle Bu, \rho v\rangle_{V^*,V} \qquad \mathrm{and} \qquad C(\rho(\cdot)y(\cdot)) = \rho(\cdot)Cy.
	\end{align*}
\end{assumption}
\noindent The most common control types satisfy this locality condition, i.e., boundary control and distributed control, as shown in the following. \changed{A similar argument also applies for distributed and boundary observations.}
\begin{lemma}[Neumann boundary control or distributed control]\label{lem:neumann_scaling}
	Let the control space $U$ and the input operator $B\in L(U,V^*)$ be defined by either of the following:
	\begin{enumerate}
		\item[i)] $U=L^2(\Gamma)$ with $\Gamma \subset \partial \Omega$ and $\langle Bu,v\rangle_{V^*,V} = \int_\Gamma u(s)\,(\operatorname{tr}v)(s) \,\mathrm{d}s$.
		\item[ii)] $U=L^2(\Omega_c)$ with $\Omega_c \subset \Omega$ and $\langle Bu,v\rangle_{V^*,V} = \int_{\Omega_c} u(\omega)\,v(\omega) \,\mathrm{d}\omega.$
	\end{enumerate} Then $B$ satisfies the condition of Assumption~\ref{as:local}.
\end{lemma}
\begin{proof}
	For case i), we compute
	\begin{align*}       \langle B(\rho u),v\rangle_{V*,V} = \int_\Gamma (\operatorname{tr}\rho)(s)\,u(s)\,(\operatorname{tr}v)(s) \,\mathrm{d}s = \int_\Gamma u(s)\,(\operatorname{tr}\rho v)(s) \,\mathrm{d}s = \langle Bu,\rho v\rangle_{V^*,V}.
	\end{align*}
	Case ii) follows analogously.
\end{proof}
\noindent We now may analyze the perturbed version of the optimality system \eqref{eq:optsys} given by
\begin{align}
\begin{pmatrix}
C^*C & 0 & A^* \\
0&I&B^*\\
A&B&0
\end{pmatrix}
\begin{pmatrix}
\tilde y\\\tilde u\\\tilde p
\end{pmatrix}
= \begin{pmatrix}
C^*Cy_d + \varepsilon_1\\\varepsilon_2\\f + \varepsilon_3
\end{pmatrix},
\end{align}
where the residual $\varepsilon = (\varepsilon_1,\varepsilon_2,\varepsilon_3)\in V^*\times U\times V^*$ satisfies 
$$\langle \varepsilon,\rho v\rangle_{V^*\times U\times V^*,V\times U\times V} = \langle \rho \varepsilon,v\rangle_{V^*\times U\times V^*,V\times U\times V} $$ 
for $\rho \in W^{1,\infty}(\Omega)$ and all $v \in V\times U\times V$. This right-hand side $\varepsilon$ may model discretization errors or perturbations of the boundary values due to, e.g.,  a domain decomposition approach. Thus, $(\tilde y,\tilde u,\tilde p)$ may be a numerical approximation of the optimal triple $(y^*,u^*,p)$ solving \eqref{eq:optsys}. Subtracting the two equations and leveraging linearity, the difference variables $(\delta y,\delta u,\delta p) = (\tilde y-y^*,\tilde u-u^*,\tilde p-p)$ satisfy
\begin{align}\label{eq:diffsys}
\begin{pmatrix}
C^*C & 0 & A^* \\
0& I&B^*\\
A&B&0
\end{pmatrix}
\begin{pmatrix}
\delta y\\\delta u\\\delta p
\end{pmatrix}
= \begin{pmatrix}
\varepsilon_1\\\varepsilon_2\\\varepsilon_3
\end{pmatrix}.
\end{align}

\noindent The following theorem provides the main result of this work: If $\varepsilon_i$, $i\in \{1,2,3\}$ are exponentially localized, then also the difference variables are exponentially localized. More precisely, we consider the following notion of exponential localization. 

\begin{definition}[Exponentially localized]\label{def:local}
	\changed{A function $f\in H^1(\mathbb{R}^d)$ is called exponentially localized around a point $z\in \mathbb{R}^d$ if there is $\mu > 0$ and a constant $c\geq 0$ such that for all domains $\Omega \subset \mathbb{R}^d$ 
		we have $\|e^{\mu \|z-\cdot\|_1}f(\cdot)\|_{H^1(\Omega)} \leq c$. Similarly, a functional $f\in H^1(\mathbb{R}^d)^*$ is called exponentially localized around a point $z\in \mathbb{R}^d$ if there is $\mu > 0$ and a constant $c\geq 0$ such that for all domains $\Omega \subset \mathbb{R}^d$ 
		we have $\|e^{\mu \|z-\cdot\|_1}f(\cdot)\|_{H^1(\Omega)^*} \leq c$.}
\end{definition}
\noindent We briefly present an interpretation of this definition. The function $f\in H^1(\mathbb{R}^d)$ is scaled with a weight that is exponentially increasing when moving away from $z\in \mathbb{R}^d$. Thus, to be bounded uniformly in $n$ for all domains $\Omega$ in the integral norm, the function $f$ has to decrease exponentially in the $1$-norm distance to $z$ with a higher rate than $\mu$. For functionals, this interpretation has to be understood in the sense that the support of the functional is exponentially localized.

\begin{theorem}\label{thm:scaling}
	Let $A$ be as in Assumption~\ref{as:secondorder} and let $B$ and $C$ satisfy Assumption~\ref{as:local}. Set $\rho(r) = e^{\mu r}$ with $\mu \in \mathbb{R}_{>0}$ such that 
	\begin{align}\label{eq:condmu}
	\mu\max\{\|a\|_{L^\infty(\Omega;\mathbb{R}^d)},\|b\|_{L^\infty(\Omega;\mathbb{R}^d)}\} + \mu^2\|a\|_{L^\infty(\Omega)} < \frac{1}{2\|M^{-1}\|}
	\end{align} and let $(\delta y,\delta u,\delta p)$ solve~\eqref{eq:diffsys}.  Assume that the perturbations of the right-hand side are exponentially localized around $z\in \Omega$ in the sense of Definition~\ref{def:local}, i.e., there is such that $\|\rho(\|z-\cdot\|_1)\varepsilon(\cdot)\|_{V^*\times U\times V^*} \leq e$ with $e\geq 0$ independent of $|\Omega|$.   Then 
	\begin{align}\label{eq:upperbound}
	\|\rho(z-\cdot)\delta y(\cdot)\|_{V} + \|\rho(z-\cdot)\delta u(\cdot)\|_{U} + \|\rho(z-\cdot)\delta p(\cdot)\|_{V} \leq \frac{\|M^{-1}\|}{2} e.
	\end{align}
	
	\noindent In particular, if the Assumptions (i) and (ii) of Theorem~\ref{thm:opnorm} are satisfied, and if $\|a\|_{L^\infty(\Omega)},\|b\|_{L^\infty(\Omega,\mathbb{R}^d)}$ are uniformly bounded in $|\Omega|$, the scaling constant $\mu>0$ and the upper bound in \eqref{eq:upperbound} can be chosen independently of $|\Omega|$ such that the perturbation variables $\delta y$, $\delta u$ and $\delta p$ are also exponentially localized around $z\in \Omega$. 
\end{theorem}
\begin{proof}
	First, we define scaled variables $\widetilde{\delta y} = \rho(z-\cdot)\delta y(\cdot)$, $\widetilde{\delta u} = \rho(z-\cdot)\delta u(\cdot)$ and $\widetilde{\delta p} = \rho(z-\cdot)\delta p(\cdot)$.
	Using Lemma~\ref{lem:perturbations} and locality of input and output maps in the sense of Assumption~\ref{as:local}, the scaled variables satisfy
	\begin{align}\label{eq:scaledvars}
	\left(
	\underbrace{\begin{pmatrix}
		C^*C & 0 & A^* \\
		0&\lambda I&B^*\\
		A&B&0
		\end{pmatrix}}_{=M} + 
	\mu 
	\underbrace{\begin{pmatrix}
		0 & 0 & F_1^* \\
		0&0&0\\
		F_1&0&0
		\end{pmatrix}}_{=:\mathcal{F}_1} + 
	\mu^2
	\underbrace{\begin{pmatrix}
		0 & 0 & F_2^* \\
		0&0&0\\
		F_2&0&0
		\end{pmatrix}}_{=:\mathcal{F}_2}
	\right)
	\begin{pmatrix}
	\widetilde {\delta y}\\\widetilde{\delta u}\\ \widetilde{\delta p}
	\end{pmatrix}
	= \begin{pmatrix}
	\tilde{\varepsilon}_1\\\tilde{\varepsilon}_2\\\tilde{\varepsilon}_3
	\end{pmatrix}
	\end{align}
	with $\tilde\varepsilon_i := \rho(z-\cdot)\varepsilon_i(\cdot)$ for $i=1,2,3$.
	We now denote $\widetilde{\delta z} = (\widetilde{\delta y},\widetilde{\delta u},\widetilde{\delta p})$ and $\tilde \varepsilon = (\tilde{\varepsilon}_1,\tilde{\varepsilon}_2,\tilde{\varepsilon}_3)$. Then, \eqref{eq:scaledvars} holds, if and only if 
	\begin{align}\label{eq:almostthere}
	\left(I + M^{-1}\left(\mu \mathcal{F}_1 + \mu^2 \mathcal{F}_2\right)\right)\widetilde{\delta z} = M^{-1}\tilde \varepsilon .
	\end{align}
	The bounds on the operator norms of $F_1$ and $F_2$ proven in~\eqref{eq:F_bdd_1norm} imply
	\begin{align*}
	\|\mathcal{F}_1\|_{L(V\times U\times V,V^*\times U\times V^*)} \leq \max\{\|a\|_{L^\infty(\Omega)},\|b\|_{L^\infty(\Omega;\mathbb{R}^d)}\}, \quad \|\mathcal{F}_2\|_{L(V\times U\times V,V^*\times U\times V^*)} \leq \|a\|_{L^\infty(\Omega)}.
	\end{align*}
	Thus, choosing $\mu\in \mathbb{R}$ such that \changed{\eqref{eq:condmu} holds} implies that
	\begin{align*}
	\changed{\|M^{-1}\left(\mu \mathcal{F}_1 + \mu^2 \mathcal{F}_2\right)\| \leq \|M^{-1}\|\left(\mu \max\{\|a\|_{L^\infty(\Omega)},\|b\|_{L^\infty(\Omega;\mathbb{R}^d)}\} + \mu^2\|a\|_{L^\infty(\Omega)}\right) < \frac{1}{2}.}
	\end{align*}
	Hence, we may apply a Neumann series argument to deduce
	\begin{align*}
	\|\left(I + M^{-1}\left(\mu \mathcal{F}_1 + \mu^2 \mathcal{F}_2\right)\right)^{-1}\| \leq \sum_{i=0}^\infty \left(\frac{1}{2}\right)^i = 2.
	\end{align*}
	Together with \eqref{eq:almostthere} this implies
	\begin{align*}
	\|\widetilde z\|_{V\times U\times V} \leq \frac{\|M^{-1}\|}{2}\|\tilde {\varepsilon}\|_{V^*\times U\times V^*} \leq \frac{\|M^{-1}\|}{2}e
	\end{align*}
	which yields the claim.
\end{proof}

\subsection{Discussion of assumptions and examples}
\noindent Before providing numerical results in the subsequent section, we illustrate the 
{stabilizability and detectability} assumption of Theorem~\ref{thm:opnorm} in view of the motivational examples of Section~\ref{sec:introduction} extended to higher space dimensions.

\begin{example}\label{ex:laplace}
	Let $V=H^1(\Omega)$ (in the case of homogeneous Neumann boundary conditions) or $V=H^1_0(\Omega)$ (in the case of homogeneous Dirichlet boundary conditions) and $\langle Ay,v\rangle_{V^*,V}:= \langle \nabla y,\nabla v\rangle_{L^2(\Omega)} + \langle cy,v\rangle_{L^2(\Omega)}$ with $c\in L^{\infty}(\Omega)$ and $\langle Bu,v\rangle_{V^*,V} = \langle u,v\rangle_{L^2(\Omega_c)}$, where $\Omega_c\subset \Omega$ is a control region.
	Then, Assumption (i) of Theorem~\ref{thm:opnorm} is satisfied if we can compensate for a nonpositive coefficient function $c$ by means of a control action, i.e., there is $r>0$ such that
	\begin{align}\label{eq:condition}
	S^{\leq r} := \{\omega\in \Omega \,|\, c(\omega) \leq r\} \subset \Omega_c.
	\end{align}
	A possible configuration is depicted in Figure~\ref{fig:config1}.
	\begin{figure}[htb]
		\begin{center}
			\begin{tikzpicture}
			\draw [thick]  plot[smooth, tension=.7] coordinates {(-4,2.5) (-3,3) (-2,2.8) (-0.8,2.2) (-0.5,1.2) (-1.5,-0) (-3.5,.5) (-4.9,1.3) (-4,2.5)};
			\draw[] plot[smooth, tension=.7] coordinates {(-4,2.5) (-2,1.5) (-2.5,0.13)};
			\draw[] plot[smooth, tension=.7] coordinates {(-3,1.5) (-2.2,1.2) (-2.5,0.6) (-3,.8) (-3.2,1.4)(-3,1.5)};
			\node at (-1.5,3.) {$\Omega$};
			\node at (-4,1.5) {$\Omega_c$};
			\node at (-2.7,1.1) {$S^{\leq r}$};
			\end{tikzpicture}
		\end{center}
		\caption{Configuration of instability region fully contained in control region.}
		\label{fig:config1}
	\end{figure}
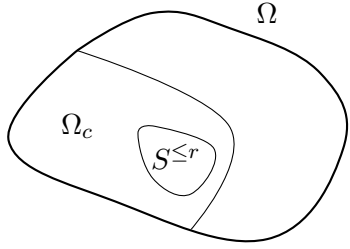
	
	\noindent More precisely, to verify Assumption~(i) of Theorem~\ref{thm:opnorm}, we may choose an arbitrary $\delta >0$ and define the multiplication operator $K_B(\omega) = \chi_{S^{\leq r}}(\omega)(-c(\omega)+\delta)$, where $\chi_{S^{\leq r}}$ is the characteristic function of $S^{\leq r}$. 
	Then,
	\begin{align*}
	\langle (A+BK_B)y,&y\rangle_{V^*,V} \\&= \langle \nabla y,\nabla y\rangle_{L^2(\Omega;\mathbb{R}^d)} + \langle cy,y\rangle_{L^2(\Omega)} + \langle \chi_{S^{\leq r}}(\cdot)\cdot(-c(\cdot)+\delta)y, y\rangle_{L^2(\Omega_c)}\\
	&= \|\nabla y\|_{L^2(\Omega;\mathbb{R}^d)}^2 + \int_{S^{\leq r}} c(\omega) y^2(\omega) + (-c(\omega) + \delta)y^2(\omega)\,\mathrm{d}\omega + \int_{\Omega \setminus S^{\leq r}} c(\omega) y^2(\omega)\,\mathrm{d}\omega\\
	&= \|\nabla y\|_{L^2(\Omega,\mathbb{R}^d)}^2 + \delta \|y\|^2_{L^2(S^{\leq r})} + r \|y\|^2_{L^2(\Omega \setminus S^{\leq r})}\\
	&\geq \min\{1,\delta,r\} \|y\|^2_{H^1(\Omega)}.
	\end{align*}
	
	
	
	
	Let us now consider a slight relaxation of \eqref{eq:condition}, i.e., assume there is $C>0$ independent of $|\Omega|$ such that
	\begin{align}\label{eq:condition2}
	\bigg| \{\omega\in \Omega \,|\, c(\omega) \leq r\} \cap \Omega \setminus \Omega_c \bigg| = \bigg|  S^{\leq r} \setminus (S^{\leq r} \cap \Omega_c)\bigg| \leq C.
	\end{align}%
	A possible geometric configuration is depicted in Figure~\ref{fig:config2}.
	\begin{figure}[htb]
		\begin{center}
			\begin{tikzpicture}[scale = 1.5]
			\draw[thick]  plot[smooth, tension=.7] coordinates {(-5.5,2.5) (-5,3) (-2,2.8) (-0.8,2.2) (-0.5,1.2) (-1.5,-0) (-3.5,.5) (-4.9,1.3) (-5.5,2.5)};
			\draw[] plot[smooth, tension=.9] coordinates {(-4,2.5) (-2,2.5) (-2.5,1) (-4.5,1.9) (-4,2.5) };
			\draw[] plot[smooth, tension=.7] coordinates {(-3,2.2) (-1.5,1.5) (-1.3,.7) (-2.5,0.5) (-3.3,.8) (-3.5,2.1)(-3,2.2)};
			\node at (-1.5,3.) {$\Omega$};
			\node at (-4,2.2) {$\Omega_c$};
			\node at (-2.8,1.6) {$S^{\leq r}\cap \Omega_c$};
			\node at (-2.3,.8) {$S^{\leq r}\setminus (S^{\leq r}\cap \Omega_c)$};
			\end{tikzpicture}
		\end{center}
		\caption{Configuration of instability region not fully contained in the control region.}
		\label{fig:config2}
	\end{figure}
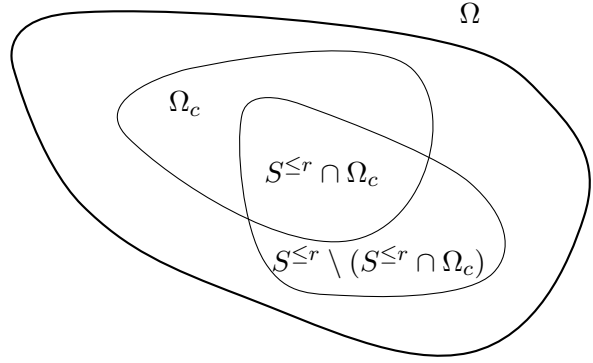%
	\noindent In this case, as in the previous example, we choose the feedback operator $K_B = \chi_{S^{\leq r}}(\omega)(-c(\omega)+\delta)$. Here, however we do not necessarily have that $S^{\leq r}\subset \Omega_c$. Thus, we compute
	\begin{align*}
	\langle (A+&BK_B)y,y\rangle_{V^*,V} \\&= \langle \nabla y,\nabla y\rangle_{L^2(\Omega;\mathbb{R}^d)} + \langle cy,y\rangle_{L^2(\Omega)} + \langle \chi_{S^{\leq r}}(-c+\delta)y, y\rangle_{L^2(\Omega_c)}\\
	&= \|\nabla y\|_{L^2(\Omega;\mathbb{R}^d)}^2 + \int_{S^{\leq r} \cap \Omega_c} c(\omega) y^2(\omega) + (-c(\omega) + \delta)y^2(\omega)\,\mathrm{d}\omega + \int_{\Omega \setminus (S^{\leq r}\cap \Omega_c)} c(\omega) y^2(\omega)\,\mathrm{d}\omega\\
	&\geq \|\nabla y\|_{L^2(\Omega;\mathbb{R}^d)}^2 +\int_{S^{\leq r} \cap \Omega_c} \delta y^2(\omega)\,\mathrm{d}\omega + r\|y\|^2_{L^2(\Omega \setminus S^{\leq r})} + \int_{S^{\leq r} \setminus (S^{\leq r}\cap \Omega_c)} c(\omega) y^2(\omega)\,\mathrm{d}\omega\\
	&\geq \min\{1,r\} \|y\|^2_{H^1(\Omega \setminus S^{\leq r})} +\|\nabla y\|_{L^2(S^{\leq r};\mathbb{R}^d)}^2 + \int_{S^{\leq r} \cap \Omega_c} \!\!\!\!\!\!\!\delta y^2(\omega)\,\mathrm{d}\omega + \int_{S^{\leq r} \setminus (S^{\leq r}\cap \Omega_c)} \!\!\!\!\!\!\!\!\!\!\!c(\omega) y^2(\omega)\,\mathrm{d}\omega.
	\end{align*}
	The first term already contains a $H^1(\Omega)$-norm on a part of the domain with a constant uniform in the domain size. Thus, we focus on the remainder and estimate
	\begin{align*}   
	&\|\nabla y\|_{L^2(S^{\leq r};\mathbb{R}^d))}^2 + \int_{S^{\leq r} \cap \Omega_c} \delta y^2(\omega)\,\mathrm{d}\omega + \int_{S^{\leq r} \setminus (S^{\leq r}\cap \Omega_c)} c(\omega) y^2(\omega)\,\mathrm{d}\omega \\&\geq 
	\!\|\nabla y\|_{L^2(S^{\leq r};\mathbb{R}^d))}^2 \!+\!\! \int_{S^{\leq r} \cap \Omega_c}\!\!\!\! \delta y^2(\omega)\,\mathrm{d}\omega - |S^{\leq r}\! \setminus\! (S^{\leq r}\cap \Omega_c)|\|c\|_{L^\infty(S^{\leq r} \setminus (S^{\leq r}\cap \Omega_c))}\|y\|^2_{L^2(S^{\leq r} \setminus (S^{\leq r}\cap \Omega_c)}\\
	&\geq \min\{C_F,\delta\}\|y\|_{H^1(S^{\leq r})}^2  - |S^{\leq r} \setminus (S^{\leq r}\cap \Omega_c)|\|c\|_{L^\infty(S^{\leq r} \setminus (S^{\leq r}\cap \Omega_c))}\|y\|_{L^2(S^\leq r)}^2\\
	&\stackrel{\eqref{eq:condition2}}{\geq} \min\{C_F - C\|c\|_{L^\infty(S^{\leq r} \setminus (S^{\leq r}\cap \Omega_c)))},\delta -C\|c\|_{L^\infty(S^{\leq r} \setminus (S^{\leq r}\cap \Omega_c))}\} \|y\|_{H^1(S^{\leq r})}^2,
	\end{align*}
	where we used the generalized Friedrichs inequality with a constant $C_F\geq 0$ depending only on $S^{\leq r} \setminus (S^{\leq r}\cap\Omega_c)$ in the second last estimate. Thus, choosing $\delta > 0$ large enough independently of $|\Omega|$ and if $C\|c\|_{L^\infty(S^{\leq r} \setminus (S^{\leq r}\cap \Omega_c))} < C_F$, we combine the above chain of inequalities to obtain Assumption (i) of Theorem~\ref{thm:opnorm}. Note that a similar argumentation can also performed by means of a Poincaré inequality, if $S^{\leq r}$ is adjacent to a boundary with homogeneous Dirichlet conditions.
\end{example}

\noindent We briefly would like to stress that deriving coercivity of the involved operators via Poincaré inequalities as in standard applications of the Lax-Milgram theorem is not admissible in view of the domain-uniformormity of the stabilizability and detectability condition in Assumption (i) of Theorem~\ref{thm:opnorm}.   
This is due to the Poincaré constant depending linearly on the diameter of the domain, cf.\ \cite{Bebe03}. As a consequence, the solution operator norm of the Poisson equation with homogeneous Dirichlet boundary conditions mapping the data to the solution measured in $L^p(\Omega)$-norms, $p<\infty$, is not uniformly bounded in the domain size, cf.\ Section~\ref{sec:introduction}.

\begin{example}[Coupled problems]
	We briefly discuss the case of elliptic problems coupled at an interface. To this end, we consider a simple example of a coupled system as depicted in Figure~\ref{fig:coupled}:
	\begin{align*}
	-\Delta y_i + \chi_{\Omega_{i,c}} u_i = f_i \quad \mathrm{on} \quad \Omega_i
	\end{align*}
	for $i\in \{1,2\}$. Here, $\chi_{\Omega_{i,c}}$ is the characteristic functions of the respective control region $\Omega_{i,c}$, $i\in \{1,2\}$.
	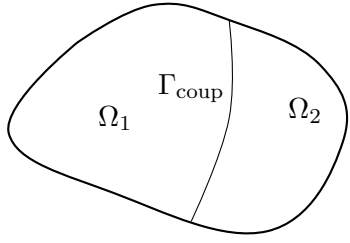
\begin{figure}[htb]
		\centering
		\begin{tikzpicture}
		\draw [thick]  plot[smooth, tension=.7] coordinates {(-4,2.5) (-3,3) (-2,2.8) (-0.8,2.2) (-0.5,1.2) (-1.5,-0) (-3.5,.5) (-4.9,1.3) (-4,2.5)};
		\draw[] plot[smooth, tension=.7] coordinates {(-2,2.8) (-2,1.5) (-2.5,0.13)};
		\node at (-2.5,1.9) {$\Gamma_\mathrm{coup}$};
		\node at (-1.,1.6) {$\Omega_2$};
		\node at (-3.5,1.5) {$\Omega_1$};
		\end{tikzpicture}
		\caption{Depiction of coupled domains.}
		\label{fig:coupled}
	\end{figure}
	
	\noindent At the coupling interface, we impose continuity and conservation of mass (or forces, depending on the application), leading to 
	\begin{align}\label{eq:coupling_cond}
	y_1 = y_2 \quad \mathrm{and}\quad \frac{\partial y_1}{\partial \nu_1} = -\frac{\partial y_2}{\partial \nu_2} \quad \mathrm{on}\ \Gamma_\mathrm{coup},
	\end{align}
	where $\nu_1$ and $\nu_2$ denote the outer unit normal of $\Omega_1$ and $\Omega_2$, respectively.
	Further, we impose homogeneous Dirichlet boundary conditions
	\begin{align*}
	y_1 = 0\quad \mathrm{on}\quad \Gamma_{1,\mathrm{D}}\subset \partial \Omega_1\setminus \Gamma_\mathrm{coup},\qquad y_2 = 0\quad \mathrm{on}\quad \Gamma_{2,\mathrm{D}}:= \partial \Omega_2  \setminus \Gamma_\mathrm{coup}.
	\end{align*}
	The weak forms of the above PDEs are given by 
	\begin{align}\label{eq:defnA1}
	\int_{\Omega_1}\nabla y_1^\top \nabla v \,\mathrm{d}\omega - \int_{\Gamma_\mathrm{coup}}\frac{\partial y_1}{\partial \nu_1}v \,\mathrm{d}s + \int_{\Omega_{1,c}} u_1 v\,\mathrm{d}\omega = \int_{\Omega_1}f_1v\,\mathrm{d}\omega \qquad \forall v\in H_{\Gamma_{1,\mathrm{D}}}^1(\Omega_1)
	\end{align}
	and
	\begin{align}\label{eq:defnA2}
	\int_{\Omega_2}\nabla y_2^\top \nabla v \,\mathrm{d}\omega + \int_{\Gamma_\mathrm{coup}}\frac{\partial y_2}{\partial \nu_1}v \,\mathrm{d}s + \int_{\Omega_{2,c}} u_2 v\,\mathrm{d}\omega = \int_{\Omega_2}f_2v\,\mathrm{d}\omega \qquad \forall v\in H_{\Gamma_{2,\mathrm{D}}}^1(\Omega_2).
	\end{align}
	The different sign in the boundary integral along $\Gamma_\mathrm{coup}$ is due to the orientation of the outer unit normals $\nu_1,\nu_2$ of $\Omega_1,\Omega_2$ at $\Gamma_\mathrm{coup}$.
	
	We now verify Assumption (i) of Theorem~\ref{thm:opnorm} for the coupled problem. In particular, this illustrates that this example may be straightforwardly extended to coupling arbitrary many domains and with constants uniform in the number of domains.
	To this end, choosing a positive but arbitrary feedback gain $k > 0$, we define multiplication-type feedback operators $K_{B_i}y_i := k y_i$ for $i\in \{1,2\}$. As a consequence, we compute, using the generalized Friedrichs inequality as in Example~\ref{ex:laplace}, the estimates
	\begin{align*}
	\langle (A_1 + \chi_{\Omega_{1,c}}K_{B_1}) y_1,y_1\rangle_{H^1_{\Gamma_{1,\mathrm{D}}}(\Omega_1)^*,H^1_{\Gamma_{1,\mathrm{D}}}(\Omega_1)} &= \|\nabla y_1\|_{L^2(\Omega_1)}^2 - \int_{\Gamma_\mathrm{coup}} \frac{\partial y_1}{\partial \nu} y_1\,\mathrm{d}s + k \|y_1\|_{L^2(\Omega_{1,c})}^2 \\
	&\geq c(|\Omega_1|,|\Omega_{1,c}|) \|y_1\|^2_{H^1(\Omega_1)} - \int_{\Gamma_\mathrm{coup}} \frac{\partial y_1}{\partial \nu_1} y_1\,\mathrm{d}s\\
	\mathrm{and}\quad \langle (A_2+ \chi_{\Omega_{2,c}}K_{B_2}) y_2,y_2\rangle_{H^1_{\Gamma_{2,\mathrm{D}}}(\Omega_2)^*,H^1_{\Gamma_{2,\mathrm{D}}}(\Omega_2)} 
	&= \|\nabla y_2\|_{L^2(\Omega_2)}^2 + \int_{\Gamma_\mathrm{coup}} \frac{\partial y_2}{\partial \nu_1} y_2\,\mathrm{d}s + k \|y_2\|_{L^2(\Omega_{2,c})}^2 
	\\&\geq c(|\Omega_2|,|\Omega_{2,c}|) \|y_2\|^2_{H^1(\Omega_2)} + \int_{\Gamma_\mathrm{coup}} \frac{\partial y_2}{\partial \nu} y_2\,\mathrm{d}s,
	\end{align*}
	where $A_1: H^1_{\Gamma_{1,\mathrm{D}}}(\Omega_1)\to H^1_{\Gamma_{1,\mathrm{D}}}(\Omega_1)^*$ and $A_2: H^1_{\Gamma_{2,\mathrm{D}}}(\Omega_2)\to H^1_{\Gamma_{2,\mathrm{D}}}(\Omega_2)^*$ are operators defined by the first two integrals in \eqref{eq:defnA1} and \eqref{eq:defnA2}, respectively.

	Adding the above inequalities and invoking the coupling conditions at $\Gamma_\mathrm{coupling}$, cf.\ \eqref{eq:coupling_cond}, we get 
	\begin{align}\label{eq:ellip1}
	\begin{split}
	\langle (A_1 + \chi_{\Omega_{1,c}}K_{B_1}) y_1,y_1\rangle_{H^1_{\Gamma_{1,\mathrm{D}}}(\Omega_1)^*,H^1_{\Gamma_{1,\mathrm{D}}}(\Omega_1)} &+ \langle (A_2+ \chi_{\Omega_{2,c}}K_{B_2}) y_2,y_2\rangle_{H^1_{\Gamma_{2,\mathrm{D}}}(\Omega_2)^*,H^1_{\Gamma_{2,\mathrm{D}}}(\Omega_2)} \\&\geq c(|\Omega_1|,|\Omega_{1,c}|) \|y_1\|^2_{H^1(\Omega_1)} + c(|\Omega_2|,|\Omega_{2,c}|) \|y_2\|^2_{H^1(\Omega_2)}.
	\end{split}
	\end{align}
	We now define \begin{align*}
	y(\omega) := \begin{cases}
	y_1(\omega) \qquad \omega\in \Omega_1\\
	y_2(\omega) \qquad \omega\in \Omega_2,
	\end{cases}
	\end{align*} 
	the Dirichlet boundary $\Gamma_{\mathrm{D}} :=\Gamma_{1,\mathrm{D}}\cup \Gamma_{2,\mathrm{D}}$ and the control region $\Omega_c :=\Omega_{1,c}\cup\Omega_{2,c}$. If we denote the Laplace operator in weak form on the coupled domain $\Omega$ by $A:H^1_{\Gamma_\mathrm{D}}(\Omega)\to H^1_{\Gamma_\mathrm{D}}(\Omega)^*$ and use \eqref{eq:ellip1},
	\begin{align*}
	\langle (A &+ \chi_{\Omega_{c}}K_{B}) y,y\rangle_{H^1_{\Gamma_\mathrm{D}}(\Omega)^*,H^1_{\Gamma_\mathrm{D}}(\Omega)} 
	\\&= \langle (-\Delta + \chi_{\Omega_{1,c}}K_{B_1}) y_1,y_1\rangle_{H^1_{\Gamma_\mathrm{D}}(\Omega_1)^*,H^1_{\Gamma_\mathrm{D}}(\Omega_1)} + \langle (-\Delta+ \chi_{\Omega_{2,c}}K_{B_2}) y_2,y_2\rangle_{H^1_{\Gamma_\mathrm{D}}(\Omega_2)^*,H^1_{\Gamma_\mathrm{D}}(\Omega_2)} 
	\\& \geq \min\{c(|\Omega_1|,|\Omega_{1,c}|),c(|\Omega_2|,|\Omega_{2,c}|)\}\left(\|y_1\|^2_{H^1(\Omega_1)} + \|y_2\|^2_{H^1(\Omega_2)}\right) \\
	&= \min\{c(|\Omega_1|,|\Omega_{1,c}|),c(|\Omega_2|,|\Omega_{2,c}|)\}\|y\|^2_{H^1(\Omega)}.
	\end{align*}
	Here, the constant in the lower bound only depends on the size of the coupled domains, and not on the total number of domains. It is clear that this argumentation carries over to an arbitrary number of domains $\Omega_1,\ldots,\Omega_N$, as long as the corresponding minimal constant $\min_{i\in N} c(|\Omega_i|,|\Omega_{i,c}|)$ resulting from the generalized Friedrichs inequality is bounded uniformly from below in $N$.
	
\end{example}
\color{black}

\section{Numerical examples: elliptic case}\label{sec:num_ell}
\noindent We perform numerical simulations on a square domain in two dimensions with side length $L>0$, i.e., $\Omega = (0,L)^2$. Further, we consider observation and control regions $\Omega_o,\Omega_c\subset \Omega$ which will be specified later and the cost functional
\begin{align*}
\min_{u\in L^2(\Omega)} \frac12\int_{\Omega_o} y^2(\omega)\,\mathrm{d}\omega  + \frac12 \int_{\Omega_c} u^2(\omega)\,\mathrm{d}\omega
\end{align*}
subject to the prototypical PDE in strong form
\begin{align}\label{eq:strongform}
\begin{split}
-\Delta y + c y &= \chi_{\Omega_c} u \qquad \mathrm{on}\ \Omega \\
y &= 0\qquad \quad \ \,\;\mathrm{on}\ \partial \Omega,
\end{split}
\end{align}
where $\chi_{\Omega_c}$ is the characteristic function of the control region $\Omega_c$ and $c\in \mathbb{R}$. Here, due to the homogeneous Dirichlet boundary conditions, we choose $V=H^1_0(\Omega)$, and get the weak formulation
\begin{align*}
Ay - Bu = 0
\end{align*}
with $A\in L(V,V^*)$ and $B\in L(U,V^*)$ defined by
\begin{align*}
\langle Ay, v\rangle_{V,V^*} = \int_\Omega \nabla y^\top \nabla v \changed{+cyv}\,\mathrm{d}\omega, \qquad \langle Bu,v\rangle = \int_{\Omega_c}u v\,\mathrm{d}\omega.
\end{align*}
In the following subsections, we will pick up the motivational examples of Section~\ref{sec:introduction} and inspect the cases $c\equiv 0$ (Poisson equation), $c\equiv -1$ (Helmholtz equation) and $c\equiv 1$ (screened Poisson equation).
We endow the state equation with a perturbation centered around $z=(L/2,L/2)^\top\in \mathbb{R}^2$ defined by
\begin{align*}
\varepsilon(\omega) = \begin{cases}
10, &\omega\in [L/2-d,L/2+d]^2\\
0, &\mathrm{otherwise}
\end{cases}
\end{align*}
for $d=2$, which clearly satisfies $\|e^{\mu \|\cdot-z\|_1}\varepsilon\|_{V^*}\leq \|e^{\mu \|\cdot-z\|_1}\varepsilon\|_{L^2(\Omega)} \leq c$ for a constant $c=c(d,\mu)\geq 0$ independent of $L$.

In the following, we will inspect the impact of this perturbation in an uncontrolled scenario and in optimal control. We will see that, if the conditions of Theorem~\ref{thm:opnorm} and Theorem~\ref{thm:scaling} are met, the optimally controlled equations will enjoy an exponential decay of this perturbation (in the sense as stated in Theorem~\ref{thm:scaling}) even though the uncontrolled equation does not.

\changed{For the numerical simulations, we set up the optimality system \eqref{eq:optsys}, discretize the involved weak formulations and operators by means of standard linear finite elements using {FEniCS}~\cite{AlnaBlech15} and solve the resulting linear equation system via the sparse linear solver \texttt{sparse.linalg.spsolve} from {SciPy} \cite{VirtGomm20}.}

\subsection{Distributed control and observation on the whole domain}
\begin{figure}[htb]
	\centering
	\includegraphics[width=.45\linewidth]{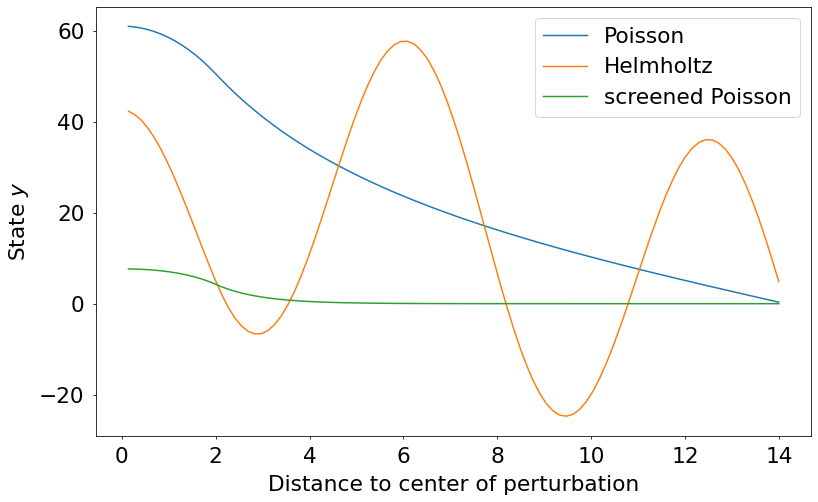}\hspace{.5cm}
	\includegraphics[width=.45\linewidth]{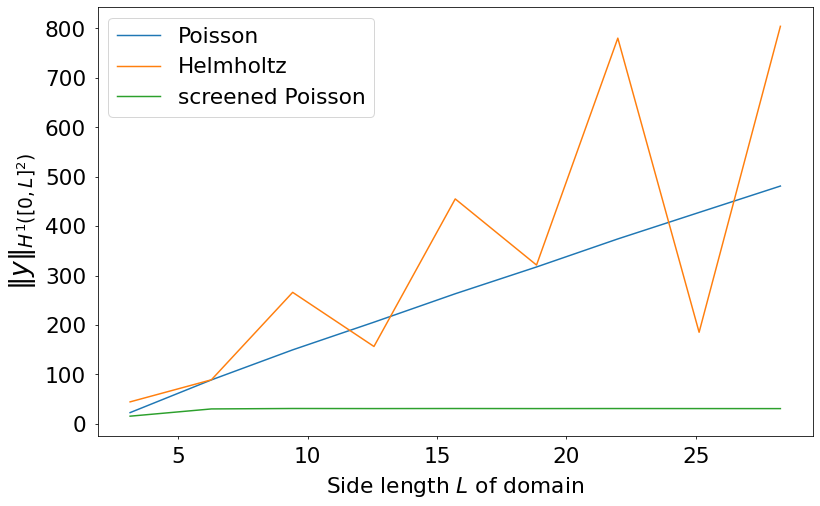}\\
	\includegraphics[width=.45\linewidth]{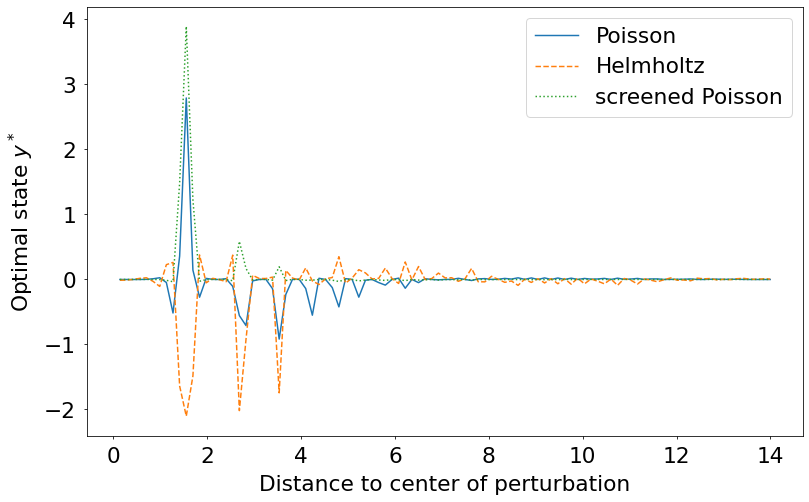}\hspace{.5cm}
	\includegraphics[width=.45\linewidth]{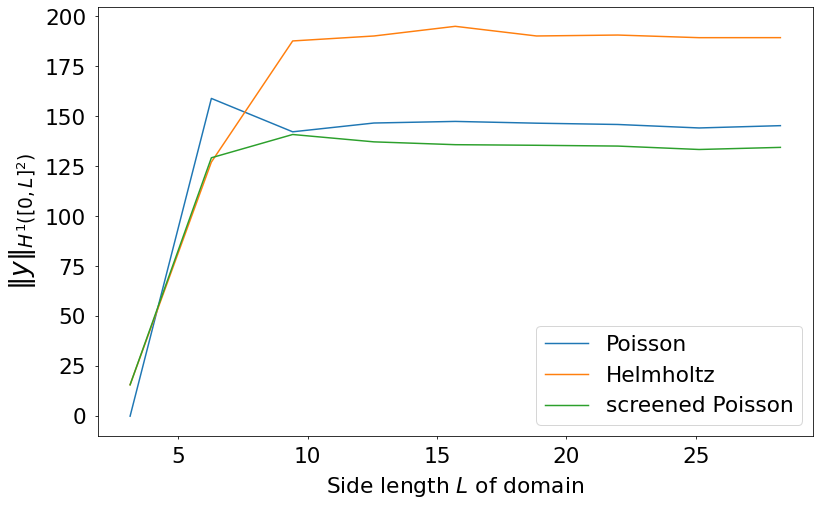}    
	\caption{Comparison of the three PDEs for a source term localized in the middle of the domain. Top: Simulation, Bottom: Optimal control.}
	\label{fig:2dres1}
\end{figure}

\noindent Let us first analyze the case $\Omega_c=\Omega_o = \Omega$, that is, control an observation on the whole domain. In view of Example~\ref{ex:laplace}, the condition \eqref{eq:condition} (and its counterpart considering the observation region) is trivially satisfied for all three equations as $\Omega \setminus \Omega_c=\Omega \setminus \Omega_o=\emptyset$. Hence we may apply Theorem~\ref{thm:opnorm} to deduce a solution operator bound uniform in the domain size, cf.\ bottom right plot of Figure~\ref{fig:2dres1} despite increasing norm of the solution of the uncontrolled problem, cf.\ top right plot of Figure~\ref{fig:2dres1}. On the left, in view of radial symmetry, we depict the uncontrolled and optimal state along a one-dimensional slice of the domain connecting $(L/2,L/2)$ and $(L/2,0)$. Due to the uniformly bounded solution operator, we may invoke the scaling result Theorem~\ref{thm:scaling} implying an exponential decay of perturbations. This is shown in the bottom left plot of Figure~\ref{fig:2dres1}.

\subsection{Distributed control on a part of the domain}
\begin{figure}[htb]
	\centering
	\includegraphics[width=.45\linewidth]{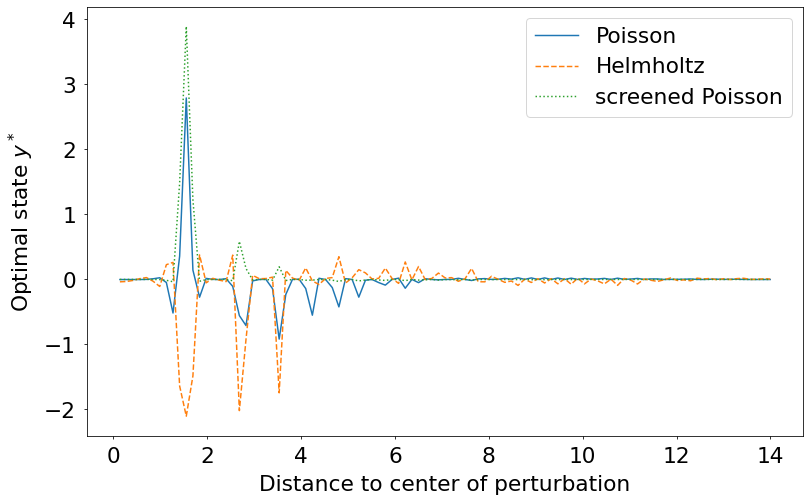}\hspace{.5cm}
	\includegraphics[width=.45\linewidth]{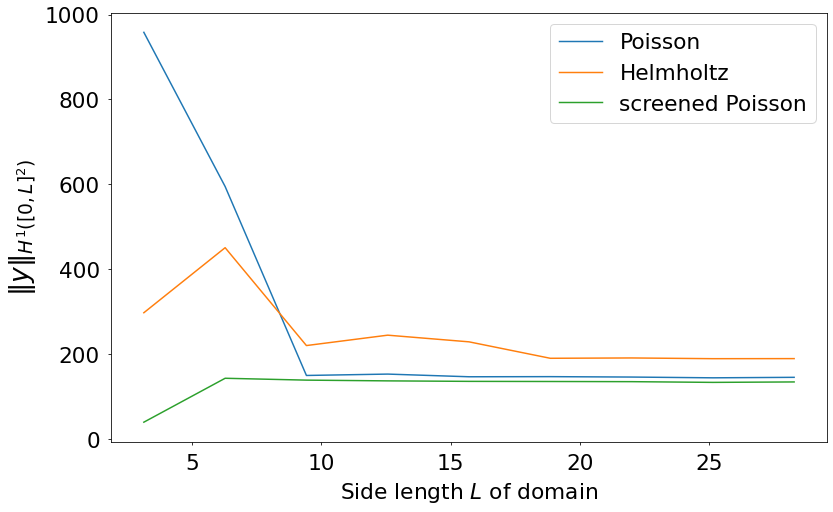}  
	\caption{Simulation of the three PDEs for a source term localized in the middle of the domain.}
	\label{fig:2dres2}
\end{figure}

\noindent Next, we analyze the case of control and observation not fully distributed over the domain. We set
\begin{align*}
\Omega_c = \Omega_o = \Omega \setminus [0,2]^2 \quad \mathrm{such\ that} \quad 
\big| \Omega \setminus \Omega_c \big| = \big| \Omega \setminus \Omega_c \big| = 4
\end{align*}
and hence, the region on which the system is not observed or controlled is bounded uniformly in the domain size. Hence, in view of Example~\ref{ex:laplace} for all three equations, conditions (i) and (ii) of Theorem~\ref{thm:opnorm} are satisfied. Thus, in Figure~\ref{fig:2dres2} we observe an exponential decay of perturbations (left) and a solution norm uniformly bounded in the domain size (right). Further, in this example, the solution norm even decreases for an increasing domain size, which is in contrast to control and observation on the full domain depicted in Figure~\ref{fig:2dres1}. We conjecture that this might stem from the increasing relative size of the control region compared to the domain size.

\subsection{Limitations}
\begin{figure}[H]
	\centering
	\includegraphics[width=.45\linewidth]{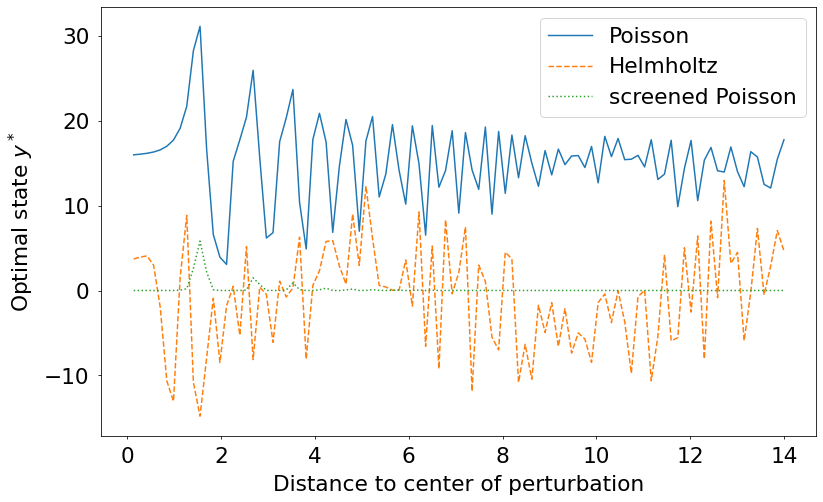}\hspace{.5cm}
	\includegraphics[width=.45\linewidth]{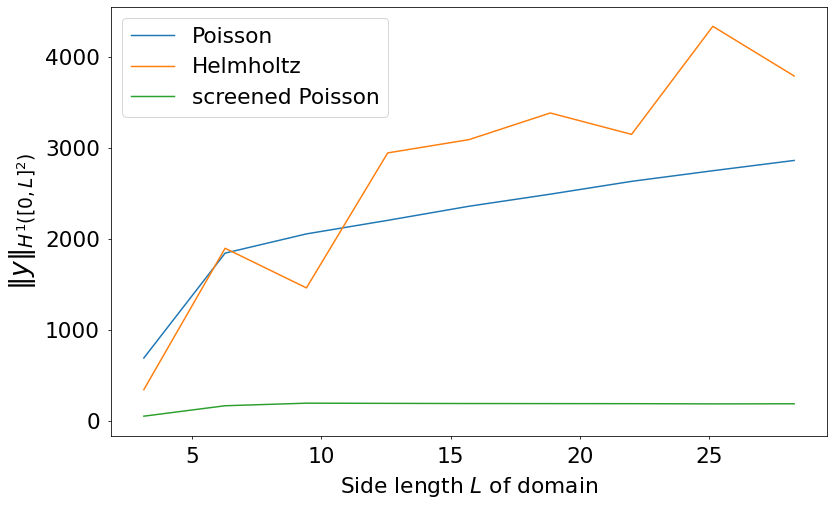}  
	\caption{Simulation of the three PDEs for a source term localized in the middle of the domain.}
	\label{fig:2dres3}
\end{figure}

\noindent Last, we show limitations on the choice of the 
input and output operators: If the unobserved or uncontrolled subdomain grows with the domain size, conditions~(i) and (ii) of Theorem~\ref{thm:opnorm} might be violated such that exponential decay can not be deduced. Here, we illustrate that this is not only a flaw in the theory, but can also be observed numerically. To illustrate this, we consider the case of control and observation on a half of the domain. More precisely, we set
\begin{align*}
\Omega_c = \Omega_o = (0,L/2)^2.
\end{align*}
In particular, for the Poisson and Helmholtz equation, conditions (i) and (ii) of Theorem~\ref{thm:opnorm} are not satisfied, as the part of the domain where the control and observation does not act is not bounded uniformly in the domain size, that is,
\begin{align*}
\big| \Omega \setminus \Omega_c \big| = \big| \Omega \setminus \Omega_o \big| = 3\left( \frac{L}{2}\right)^2 = \frac{3L^2}{4}.
\end{align*}
Thus we see that the perturbations are only exponentially damped in case of the screened Poisson equation,  cf.\ the left plot of Figure~\ref{fig:2dres3}; in this case, conditions (i) and (ii) of Theorem~\ref{thm:opnorm} are satisfied with $K_C=K_B=0$. Further, again only for the screened Poisson equation, the solution norm is bounded uniformly in the domain size, as depicted in the right plot of Figure~\ref{fig:2dres3}.

\section{Spatial decay of perturbations: Parabolic optimal control}\label{sec:parab}
\noindent In this part, we provide an extension of the localized sensitivity results proven in Theorem~\ref{thm:scaling} to parabolic problems.
To this end, we again consider $V\subset H^1(\Omega)$ endowed with the $H^1(\Omega)$-topology such that $V\hookrightarrow L^2(\Omega) \hookrightarrow V^*$ is a Gelfand triple with continuous and dense embeddings. Let $T>0$ be a time horizon and consider the dynamic OCP
\begin{align}\label{eq:OCPp}
\begin{split}
\min_{(y,u)\in W(0,T)\times L^2(0,T;U)} \frac12 (\|C(y-y_d)\|_{L^2(0,T;Y)}^2 &+ \|u\|_{L^2(0,T;U)}^2)\\
\text{s.t. } \dot y + Ay + Bu &= f \qquad \text{on } [0,T]\times \Omega\\
y(0)&= y^0,
\end{split}
\end{align}
where $B\in L(L^2(0,T;U),L^2(0,T;V^*))$, $C\in L(L^2(0,T;V),L^2(0,T;Y))$, $y^0\in L^2(\Omega)$, $y_d\in L^2(0,T;V)$, $f\in L^2(0,T;V^*)$ and $A\in L(L^2(0,T;V),L^2(0,T;V^*))$ is again a second order elliptic differential operator, i.e.,
$$
\langle Ay,v\rangle_{L^2(0,T;V^*),L^2(0,T,V)} = \langle a\nabla y,\nabla v\rangle_{L^2((0,T)\times \Omega)} + \langle b^\top\nabla y,v \rangle_{L^2((0,T)\times \Omega)} + \langle cy,v\rangle_{L^2((0,T)\times \Omega)}
$$ 
with coefficients $a\in L^\infty((0,T)\times \Omega;\mathbb{R}^{d\times d})$ s.t.~there is $\underline{a}>0$ such that $v^\top a(t,\omega)v\geq \underline{a}\|v\|^2$ for a.e.~$(t,\omega)\in (0,T)\times \Omega$ and $v\in \mathbb{R}^d$, $c\in L^\infty((0,T)\times \Omega)$ and $b\in L^\infty((0,T)\times \Omega;\mathbb{R}^d)$. It may be straightforwardly shown \cite[§6.2.2.]{evans2022partial} that $A$ satisfies the G\aa rding (or generalized ellipticity) inequality
\begin{align*}
\exists \beta>0\, \gamma\geq 0:\qquad  \langle Av,v\rangle_{L^2(0,T;V^*),L^2(0,T;V)} + \gamma\|v\|_{L^2(0,T;L^2(\Omega))}^2 \geq \beta \|v\|_{L^2(0,T,V)}^2,
\end{align*}
Understanding the time derivative in the distributional sense, we define the usual space
\begin{align*}
W(0,T) := \{v\in L^2(0,T;V)\,|\,\dot v \in L^2(0,T;V^*)\}
\end{align*}
which forms a Banach space when endowed with the norm $\|v\|_{W(0,T)} := \|v\|_{L^2(0,T;V)} + \|\dot v\|_{L^2(0,T;V^*)}$. Thus, the dynamics in \eqref{eq:OCPp} are understood in a weak sense:
\begin{align}\label{eq:dynamics}
\langle \dot y + Ay + Bu,v\rangle_{L^2(0,T;V^*),L^2(0,T;V)} = \langle f,v\rangle_{L^2(0,T;V^*),L^2(0,T;V)} \qquad \forall v\in L^2(0,T;V).
\end{align}
For each $f\in L^2(0,T;V^*)$ and $y\in L^2(\Omega)$, there is a unique weak solution of the dynamics \eqref{eq:dynamics}, cf.\ \cite{Wlok87} or \cite[Chapter 23]{Zeid13}. In view of maximal parabolic regularity, it can be shown that $\frac{\mathrm{d}}{\mathrm{d}t}+A$ together with the initial condition yields an isomorphism from $W(0,T)$ to $L^2(0,T;V^*)$. Further, we have the continuous embedding $W(0,T)\hookrightarrow C(0,T;L^2(\Omega))$ and the embedding constant does not depend on~$\Omega$ and $T$, cf.\ \cite[Problem 23.10 p.446]{Zeid13}. Further, we will make use of the integration-by-parts formula: For all $y\in W(0,T)$,
\begin{align}\label{eq:intparts}
\langle \dot y,y\rangle_{L^2(0,T;V^*),L^2(0,T;V)} = \frac{1}{2}\left(\|y(T)\|_{L^2(\Omega)}^2 - \|y(0)\|^2_{L^2(\Omega)}\right).
\end{align}
We briefly recap optimality conditions as provided in \cite{Schi13}. To this end, for $t\in [0,T]$, we define a time evaluation operator $E_t:W(0,T)\to L^2(\Omega)$ by $E_t y:= y(t)$ which is continuous and bounded due to $W(0,T)\hookrightarrow C(0,T;L^2(\Omega))$.
\begin{lemma}\label{lem:optcond_p}
	Let $(y^*,u^*)\in W(0,T)\times L^2(0,T;U)$ be optimal for \eqref{eq:OCPp}. Then there is an adjoint state $p\in W(0,T)$ such that
	\begin{align}\label{eq:optsys_p}
	\underbrace{\begin{pmatrix}
		C^*C & 0 & -\nicefrac{\mathrm{d}}{\mathrm{d}t}+A^* \\
		0&0&E_T\\
		0&I&B^*\\
		\nicefrac{\mathrm{d}}{\mathrm{d}t}+A&B&0\\
		E_0&0&0
		\end{pmatrix}}_{=:M}
	\begin{pmatrix}
	y^*\\u^*\\p
	\end{pmatrix}
	= \begin{pmatrix}
	C^*Cy_d\\0\\0\\f\\y^0
	\end{pmatrix}.
	\end{align}
\end{lemma}
\noindent Again, our main goal is to prove locality of perturbations of the optimality system. To this end, analogously to~\eqref{eq:diffsys}, we consider perturbation variables $(\delta y,\delta u,\delta p)\in W(0,T)\times L^2(0,T;U)\times W(0,T)$ solving
\begin{align}\label{eq:diffsys_p}
M
\begin{pmatrix}
\delta y\\\delta u\\\delta p
\end{pmatrix}
= \varepsilon
\end{align}
with $\varepsilon = (\varepsilon_1,\varepsilon_2,\varepsilon_3,\varepsilon_4,\varepsilon_5) \in L^2(0,T;V^*)\times L^2(\Omega) \times L^2(0,T;U)\times L^2(0,T,V^*)\times L^2(\Omega)$ satisfying $\langle \varepsilon,\rho v\rangle = \langle \rho \varepsilon,v\rangle$ for $\rho \in L^\infty(0,T;W^{1,\infty}(\Omega))$ and all $v\in L^2(0,T;V)\times L^2(\Omega)\times L^2(0,T;U)\times L^2(0,T;V)\times L^2(\Omega)$.



\subsection{Spatially exponential decay in the parabolic case}
\noindent We now provide the surrogate of Theorem~\ref{thm:opnorm} for the parabolic optimal control problem \eqref{eq:OCPp}. In this context, we will consider bounds on the solution operator that are uniform in both time horizon $T$ and size of the spatial domain $|\Omega|$. \changed{For a motivation why this uniformity in the time horizon is crucial, we refer to Subsection~\ref{subsec:crucial}.}

\begin{theorem}\label{thm:opnorm_p}
	Let the operators $A,B$ and $C$ in \eqref{eq:OCPp} be bounded uniformly in $| \Omega |$ and $T$. Suppose the existence of a constant $\alpha > 0$ independent of $|\Omega|$ and $T$ and feedback operators $K_B\in L(L^2(0,T;V),L^2(0,T;U))$ and $K_C\in L(L^2(0,T;Y),L^2(0,T;V^*))$, bounded uniformly in $|\Omega|$ and $T$ such that
	\begin{enumerate}
		\item[(i)] $(A+BK_B)y,y\rangle_{L^2(0,T,V^*),L^2(0,T;V)} \geq \alpha \|y\|_{L^2(0,T;V)}^2$ and
		\item[(ii)] $\langle (A+K_CC)y,y\rangle_{L^2(0,T;V^*), L^2(0;T,V)} \geq \alpha \|y\|_{L^2(0,T;V)}^2$
	\end{enumerate}
	hold. Then there is a constant $c\geq 0$ independent of $|\Omega|$ and $T$ such that
	\begin{align*}
	\|M^{-1}\|_{L(L^2(0,T;V^*)\times L^2(\Omega)\times L^2(0,T;U)\times L^2(0,T;V^*)\times L^2(\Omega) , W(0,T)\times L^2(0,T;U)\times W(0,T))} &\leq c
	\end{align*}
	for the operator matrix $M$ defined in the optimality system~\eqref{eq:optsys_p}.
\end{theorem}
\begin{proof}
	We proceed similarly to the proof of Theorem~\ref{thm:opnorm} and abbreviate $ (y,v) := \langle y,v\rangle_{L^2(0,T;V^*), L^2(0,T;V)}$. We omit all subscripts for operator norms as they will be clear from context. 
	Let $(y,u,p)\in W(0,T)\times L^2(0,T;U)\times W(0,T)$ solve \eqref{eq:diffsys_p}. We eliminate the control by means of the third equation of \eqref{eq:diffsys_p} via
	\begin{align}\label{eq:subsU_p}
	u = -B^* p + \varepsilon_3
	\end{align}
	and obtain the reduced system
	\begin{align}\label{eq:optcond_elim_p}\begin{pmatrix}
	C^*C & -\nicefrac{\mathrm{d}}{\mathrm{d}t}+A^*\\
	0&E_T\\
	\nicefrac{\mathrm{d}}{\mathrm{d}t} + A & -BB^*\\
	E_0&0
	\end{pmatrix}
	\begin{pmatrix}
	y\\p
	\end{pmatrix}
	=
	\begin{pmatrix}
	\varepsilon_1\\
	\varepsilon_2\\
	\varepsilon_4 - B \varepsilon_3\\
	\varepsilon_5
	\end{pmatrix}.
	\end{align}
	Using Assumption~(ii), we add $(BK_By,v)$ on both sides of the second last equation and get
	\begin{align*}
	\langle \dot y,v\rangle + ((A+K_CC)y,v) = (\varepsilon_4-B\varepsilon_3 + BB^*p + K_CCy,v).
	\end{align*}
	for all $v\in L^2(0,T;V)$. Thus, testing this equation with $y$, invoking Assumption~(i) and \eqref{eq:intparts}, we get
	\begin{align*}
	&\frac12\|y(T)\|_{L^2(\Omega)}^2 + \alpha \|y\|_{L^2(0,T;V)}^2 \\&\leq (\varepsilon_4-B\varepsilon_3 + BB^*p + K_CCy,y) + \frac12\|y(0)\|^2_{L^2(\Omega)}\\
	&\leq c\left(\|\varepsilon_3\|_{L^2(0,T;U)} + \|\varepsilon_4\|_{L^2(0,T;V^*)} + \|B^*p\|_{L^2(0,T;U)} + \|Cy\|_{L^2(0,T;Y)}\right)\|y\|_{L^2(0,T;V)} +  \frac12\|\varepsilon_5\|^2_{L^2(\Omega)}
	\end{align*}
	for a constant $c=c(\|B\|,\|K_C\|)\geq 0$. Using the simple estimate $(a+b)^2\leq 2(a^2+b^2)$ for $a,b\in \mathbb{R}$, there is a constant $c\geq 0$ independent of $|\Omega|$ such that
	\begin{align}\label{eq:p_state}
	\begin{split}
	\|&y(T)\|^2_{L^2(\Omega)} + \|y\|_{L^2(0,T;V)}^2 \\&\leq c\left(\|\varepsilon_3\|^2_{L^2(0,T;U)} + \|\varepsilon_4\|^2_{L^2(0,T;V^*)} + \|B^*p\|^2_{L^2(0,T;U)} + \|Cy\|^2_{L^2(0,T;Y)} + \|\varepsilon_5\|^2_{L^2(\Omega)} \right).
	\end{split}
	\end{align}
	
	The bound for the adjoint equation follows analogously by using Assumption (i) and we get
	\begin{align}\label{eq:p_adjoint}
	\begin{split}
	\|p(0)\|^2_{L^2(\Omega)} &+ \|p\|^2_{L^2(0,T;V)} \\&\leq c\left(\|\varepsilon_1\|^2_{L^2(0,T;V^*)} + \|B^*p\|^2_{L^2(0,T;U)} + \|Cy\|^2_{L^2(0,T;Y)} + \|\varepsilon_2\|^2_{L^2(\Omega)}  \right)
	\end{split}
	\end{align}
	with a constant $c\geq 0$ independent of $|\Omega|$.
	To eliminate the term $\|B^*p\|^2_{L^2(0,T;U)} + \|Cy\|^2_{L^2(0,T;Y)}$, we test the state equation with the adjoint state, and the adjoint equation with the state and subtract the former from the latter. This yields
	\begin{align*}
	\|Cy\|^2_{L^2(0,T;Y)} &+ \|B^*p\|^2_{L^2(0,T;Y)}\\
	&= (\tfrac{\mathrm{d}}{\mathrm{d}t} p,y) - (\tfrac{\mathrm{d}}{\mathrm{d}t} y,p) + (\varepsilon _1,y) - (\varepsilon_4-Be_3,p)\\
	&= \langle p(T),y(T)\rangle_{L^2(\Omega)}-\langle p(0),y(0)\rangle_{L^2(\Omega)} + (\varepsilon _1,y) - (\varepsilon_4-B\varepsilon_3,p)\\
	&\leq \|\varepsilon_2\|_{L^2(\Omega)}\|y(T)\|_{L^2(\Omega)} + \|\varepsilon_5\|_{L^2(\Omega)}\|p(0)\|_{L^2(\Omega)} + \|\varepsilon_1\|_{L^2(0,T;V^*)}\|y\|_{L^2(0,T;V)} \\
	&\qquad + \left(\|\varepsilon_4\|_{L^2(0,T;V^*)} + \|B\|\|\varepsilon_3\|_{L^2(0,T;U)}\right)\|p\|_{L^2(0,T;V)}.
	\end{align*}
	Combining this with \eqref{eq:p_state} and \eqref{eq:p_adjoint} and again using $(a+b)^2\leq 2(a^2+b^2)$ for $a,b\in \mathbb{R}$, there is a constant $c\geq 0$ such that
	\begin{align}\label{eq:p_combined}
	\begin{split}
	\|y&\|^2_{L^2(0,T;V)} + \|p\||^2_{L^2(0,T;V)} \\
	&\leq c \left(\|\varepsilon_1\|^2_{L^2(0,T;V^*)} + \|\varepsilon_2\|^2_{L^2(\Omega)} + \|\varepsilon_3\|^2_{L^2(0,T;U)} + \|\varepsilon_4\|^2_{L^2(0,T;V^*)} + \|\varepsilon_5\|^2_{L^2(\Omega)}     \right).
	\end{split}
	\end{align}
	Using the state and the adjoint equation, we deduce
	\begin{align*}
	\|\dot y\|_{L^2(0,T;V^*)} \leq \|A\|\|y\|_{L^2(0,T;V)}  + \|BB^*\|\|p\|_{L^2(0,T;V)}  + \|\varepsilon_4\|_{L^2(0,T;V^*)} \changed{+\|B\|\|\varepsilon_3\|_{L^2(0,T;U)}}
	\end{align*}
	and 
	\begin{align*}
	\|\dot p\|_{L^2(0,T;V^*)} \leq \|A\|\|p\|_{L^2(0,T;V)}  + \|C^*C\|\|y\|_{L^2(0,T;V)}  + \|\varepsilon_1\|_{L^2(0,T;V^*)}.
	\end{align*}
	Thus, as $\|v\|_{W(0,T)}=\|v\|_{L^2(0,T;V)} + \|\dot{v}\|_{L^2(0,T;V^*)}$, combining this with \eqref{eq:p_combined}, yields the desired $W(0,T)$-bound on the state and adjoint. The bound for the control follows by resubstitution via \eqref{eq:subsU_p}.
\end{proof}

We briefly comment on the assumptions of above theorem.
\begin{remark}
	Assumptions (i) and (ii) of Theorem~\ref{thm:opnorm_p} are sufficient to deduce exponential stabilizability of $(A,B)$ and detectability of $(A,C)$ without overshoot, that is, weak solutions to $\dot y = A+BK_By$, $y(0)=y^0$ satisfy $\|y(t)\|_{L^2(\Omega)} \leq e^{-\alpha t}\|y^0\|_{L^2(\Omega)}$, cf.\ \cite{GrunScha19} and \cite[Remark 3.19]{Scha21}.
\end{remark}

We may now deduce the main result of this section, implying spatially exponential decay in parabolic optimal control problems.

\begin{theorem}\label{thm:scaling_p}
	Assume that $B$ and $C$ act local in space in the sense of Assumption~\ref{as:local}. Further, set $\rho(r) = e^{\mu r}$ with $\mu \in \mathbb{R}_{>0}$ satisfying 
	$$
	\mu\max\{\|a\|_{L^\infty((0,T)\times\Omega);\mathbb{R}^{d\times d})},\|b\|_{L^\infty(0,T)\times\Omega;\mathbb{R}^d)}\} + \mu^2\|a\|_{L^\infty((0,T)\times\Omega;\mathbb{R}^{d\times d}))}  < \frac{1}{2\|M^{-1}\|}
	$$ 
	and let $(\delta y,\delta u,\delta p)$ solve \eqref{eq:diffsys_p}. Assume that the perturbations of the right-hand side are exponentially localized around $z\in \Omega$, 
	i.e., there is $z\in \Omega$ such that $$\|\rho(\|z-\cdot\|_1)\varepsilon\|_{L^2(0,T;V^*)\times L^2(\Omega)\times L^2(0,T;U)\times L^2(0,T;V^*)\times L^2(\Omega)} \leq e$$ with $e\geq 0$ independent of $|\Omega|$. Then 
	\begin{align}\label{eq:upperbound_p}
	\|\rho(z-\cdot)\delta y(\cdot)\|_{W(0,T)} + \|\rho(z-\cdot)\delta u(\cdot)\|_{L^2(0,T;U)} + \|\rho(z-\cdot)\delta p(\cdot)\|_{W(0,T)} \leq \frac{\|M^{-1}\|}{2} e.
	\end{align}
	
	\noindent In particular, if the Assumptions (i) and (ii) of Theorem~\ref{thm:opnorm_p} are satisfied and if $\|a\|_{L^\infty(0,T)\times\Omega;\mathbb{R}^{d\times d}))}$ and  $\|b\|_{L^\infty(0,T)\times\Omega;\mathbb{R}^d)}$ are bounded uniformly in $|\Omega|$ and $T$, the scaling constant $\mu$ and the upper bound in \eqref{eq:upperbound_p} can be chosen independently of $|\Omega|$ and $T$ such that the perturbation variables $\delta y$, $\delta u$ and $\delta p$ are also exponentially localized around $z\in \Omega$. 
\end{theorem}
\begin{proof}
	The proof is analogous to the proof of Theorem~\ref{thm:scaling}.
\end{proof}

\begin{remark}
	The result of Theorem~\ref{thm:scaling_p} also holds true when replacing the fixed point $z\in \Omega$ by a possibly time-dependent point $z:[0,T]\to \Omega$.
\end{remark}

\subsection{The necessity for a spatial decay uniform in the time horizon}\label{subsec:crucial}
\noindent In the main result presented in the previous subsection in Theorem~\ref{thm:scaling_p}, we required that the decay parameter $\mu$ entering the scaling function $\rho(r) = e^{\mu r}$ is uniform in the time horizon of the OCP. Correspondingly, in Theorem~\ref{thm:opnorm_p}, we formulated sufficient conditions to also render the solution operator norm independent of $T$. Loosely speaking, this uniformity  ensures that the spatial decay does not deteriorate for larger and larger time horizons.

To illustrate that this uniformity in the time horizon is crucial when considering locality in space, we consider the example of the heat equation with homogeneous Neumann boundary conditions. Setting $V=H^1(\Omega)$, its weak formulation is given by
\begin{subequations}\label{eq:parab_PDE}
	\begin{align}
	\label{eq:pdyn}
	\begin{split}
	\langle \dot y + Ay, v &\rangle_{L^2(0,T;H^1(\Omega)^*),L^2(0,T;H^1(\Omega))} \\&= \langle f,v\rangle_{L^2(0,T;H^1(\Omega)^*,L^2(0,T;H^1(\Omega))} \qquad \forall v\in L^2(0,T;H^1(\Omega)),
	\end{split}\\
	y(0)&=y^0,\label{eq:pinit}
	\end{align}
\end{subequations}
with $y^0\in L^2(\Omega)$ and $f\in L^2(0,T;H^1(\Omega)^*)$ where
\begin{align*}
\langle Ay,v\rangle_{L^2(0,T;H^1(\Omega)^*),L^2(0,T;H^1(\Omega))} = \langle \nabla y,\nabla v\rangle_{L^2([0,T]\times \Omega)} \qquad \forall y,v\in L^2(0,T;H^1(\Omega)).
\end{align*}
Whereas it is clear that, e.g., an initial perturbation has a global impact due to the homogeneous Neumann boundary conditions and diffusion, which loosely speaking only distributes initial temperature distribution over the domain, we will show that the solution operator of the PDE can still be bounded uniformly in the size of the spatial domain, i.e., the operator
\begin{align*}
\mathcal{A}: W(0,T) &\to L^2(0,T;H^1(\Omega)^*) \times L^2(\Omega)\\
y &\mapsto \begin{pmatrix}
\frac{\mathrm{d}}{\mathrm{d}t}y+ Ay\\
y(0)
\end{pmatrix}
\end{align*}
satisfies $\|\mathcal{A}^{-1}\|_{L(L^2(0,T;H^1(\Omega)^*) \times L^2(\Omega),W(0,T))}\leq c$ for a constant $c\geq 0$ independent of $|\Omega|$. To this end, we first introduce a scaled variable $\tilde y(t,x) = e^{-t}y(t,x)$ which, due to the product rule, satisfies
\begin{align}\label{eq:parab_subs}
\dot {\tilde y} = -\tilde y + e^{-t}\dot y
\end{align}
such that, setting $\tilde f(t,x) :=e^{-t}f(t,x)$, the dynamics \eqref{eq:parab_PDE} are equivalent to
\begin{align}\label{eq:pdynscaled}
\dot {\tilde y} + (A+ I)\tilde y = \tilde f .
\end{align}
Now we test this equation with $\tilde y$, invoke \eqref{eq:intparts}, use $\tilde y (0)=y^0$ and obtain
\begin{align*}
\frac12\left(\|\tilde y(T)\|_{L^2(\Omega)}^2 - \|y^0\|_{L^2(\Omega)}^2\right) + \|\nabla \tilde y\|_{L^2(0,T;L^2(\Omega))}^2 &+ \|\tilde y\|^2_{L^2(0,T;L^2(\Omega))} \\
&= \langle \tilde f,\tilde y\rangle_{L^2(0,T;H^1(\Omega)^*),L^2(0,T;H^1(\Omega))}
\end{align*}
which implies
\begin{align*}
\|\tilde y \|_{L^2(0,T;H^1(\Omega))} \leq \frac12 \|y^0\|_{L^2(\Omega)} + \|\tilde f\|_{L^2(0,T;H^1(\Omega)^*)}.
\end{align*}
Using \eqref{eq:pdynscaled}, we further obtain
\begin{align*}
\|\tilde y\|_{W(0,T)} &= \|\tilde y \|_{L^2(0,T;H^1(\Omega))} +  \|\dot{\tilde y} \|_{L^2(0,T;H^1(\Omega)^*)} 
\\&\leq \|\tilde y \|_{L^2(0,T;H^1(\Omega))} +  \|A+I\| \|\tilde y\|_{L^2(0,T;H^1(\Omega))} + \|\tilde f\|_{L^2(0,T;H^1(\Omega)^*)}
\\& \leq 2\|\tilde y \|_{L^2(0,T;H^1(\Omega))} + \|\tilde f\|_{L^2(0,T;H^1(\Omega)^*)}
\\& \leq \|y^0\|_{L^2(\Omega)} + 3\|\tilde f\|_{L^2(0,T;H^1(\Omega)^*)}
\end{align*}
as
\begin{align*}     
\|A+&I\|_{L(L^2(0,T;H^1(\Omega)),L^2(0,T;H^1(\Omega)^*)} \\&= \sup_{\|v\|_{L^2(0,T;H^1(\Omega))}=1}\sup_{\|y\|_{L^2(0,T;H^1(\Omega))}=1} \langle \nabla y,\nabla v\rangle + \langle y,v\rangle  
\\&\leq \sup_{\|v\|_{L^2(0,T;H^1(\Omega))}=1}\sup_{\|y\|_{L^2(0,T;H^1(\Omega))}=1}\|\nabla y\|\|\nabla v\| + \|y\|\|v\|
\\&\leq \sup_{\|v\|_{L^2(0,T;H^1(\Omega))}=1}\sup_{\|y\|_{L^2(0,T;H^1(\Omega))}=1}\frac12 \|y\|^2_{L^2(0,T;H^1(\Omega))} + \frac12 \|v\|^2_{L^2(0,T;H^1(\Omega))}\\
&= 1.
\end{align*}
Hence, resubstituting the original state we get
\begin{align*}
\|y\|_{L^2(0,T;H^1(\Omega))} = \|e^{t} \tilde y\|_{L^2(0,T;H^1(\Omega))}  \leq e^T \|\tilde y\|_{L^2(0,T;H^1(\Omega))} \leq e^T \left(\frac12 \|y^0\|_{L^2(\Omega)} + \|\tilde f\|_{L^2(0,T;H^1(\Omega)^*)}\right)
\end{align*}
and via \eqref{eq:parab_subs}, we compute
\begin{align*}
\|\dot y\|_{L^2(0,T;H^1(\Omega)^*)} = \|e^t (\tilde y + \dot{\tilde y})\|_{L^2(0,T;H^1(\Omega)^*)} \leq e^T \|\tilde y\|_{W(0,T)} \leq e^T \left(\|y^0\|_{L^2(\Omega)} + 3\|\tilde f\|_{L^2(0,T;H^1(\Omega)^*)}\right) 
\end{align*}
such that 
\begin{align*}
\|y\|_{W(0,T)} \leq e^T \left(\frac32 \|y^0\|_{L^2(\Omega)} + 4\|\tilde f\|_{L^2(0,T;H^1(\Omega)^*)}\right) \leq e^T\left(\frac32 \|y^0\|_{L^2(\Omega)} + 4\| f\|_{L^2(0,T;H^1(\Omega)^*)}\right).
\end{align*}
Thus, we have the following bound on the solution operator corresponding to the parabolic PDE \eqref{eq:parab_PDE}
\begin{align*}
\|\mathcal{A}^{-1}\|_{L(L^2(0,T;H^1(\Omega)^*) \times L^2(\Omega),W(0,T))} \leq 4e^T.
\end{align*}
As this bound grows exponentially in the time horizon, this means that, when generalizing Theorem~\ref{thm:scaling} to the parabolic case, the parameter $\mu\sim\frac{1}{\|\mathcal{A}^{-1}\|}$ corresponding to the spatial decay vanishes for $T\to \infty$. Intuitively, this corresponds to the diffusion with homogeneous Neumann boundary conditions which distributes the initial perturbation over the domain and hence in particular does not exhibit spatial locality.
\color{black}
\section{Numerical examples: parabolic case}\label{sec:num_parab}
\noindent We present numerical results for a heat equation on a one-dimensional domain $\Omega = (0,L)$, $L>0$, with distributed control and observation and homogeneous Neumann boundary conditions given by
\begin{align*}
\dot y - \Delta y + cy &= u + f\qquad \mathrm{on}\ [0,T]\times \Omega\\
\partial_\nu y &= 0 \qquad \quad \ \ \; \mathrm{on}\ [0,T]\times \partial \Omega\\
y(0)&= 0 \qquad \quad \ \ \;\mathrm{on}\ \Omega
\end{align*}
for a coefficient $c \in \{-1,0,1\}$. The one-dimensional spatial domain is chosen to provide illustrations of the results in space and time. In particular, we stress that the main result of Theorem~\ref{thm:scaling_p} holds for general spatial dimensions.

\changed{Similar to the numerical implementation for the elliptic case described in Section~\ref{sec:num_ell}, we set up the optimality system, discretize the involved operators and solve the linear equations system. Space discretization is performed via linear finite elements using {FEniCS}~\cite{AlnaBlech15} and for time discretization, we use an implicit Euler method.}

\begin{figure}[H]
	\centering
	\includegraphics[width=.32\linewidth]{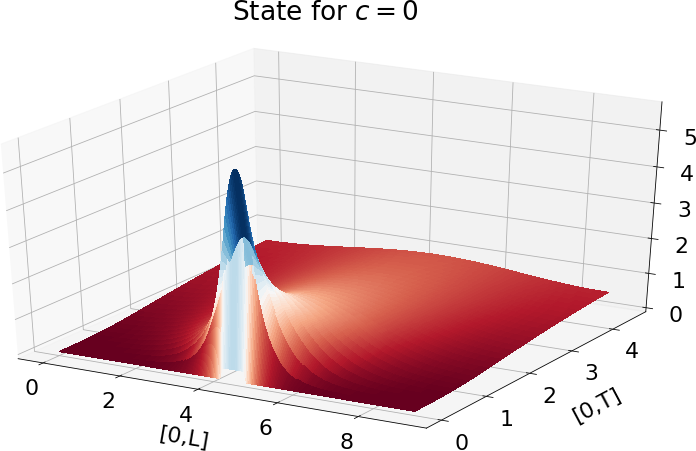}
	\includegraphics[width=.32\linewidth]{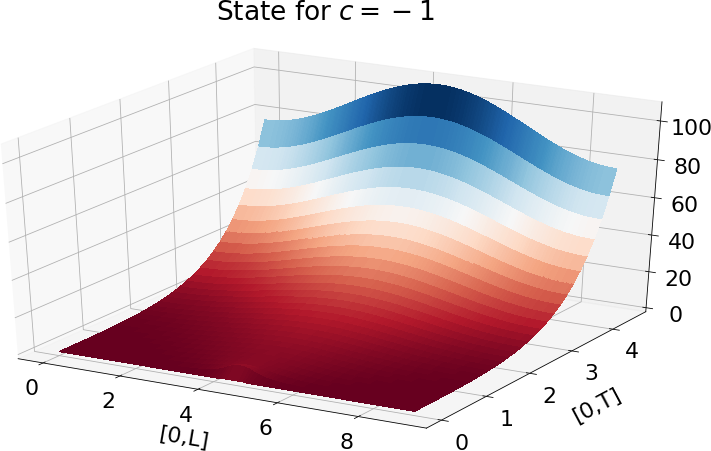}
	\includegraphics[width=.32\linewidth]{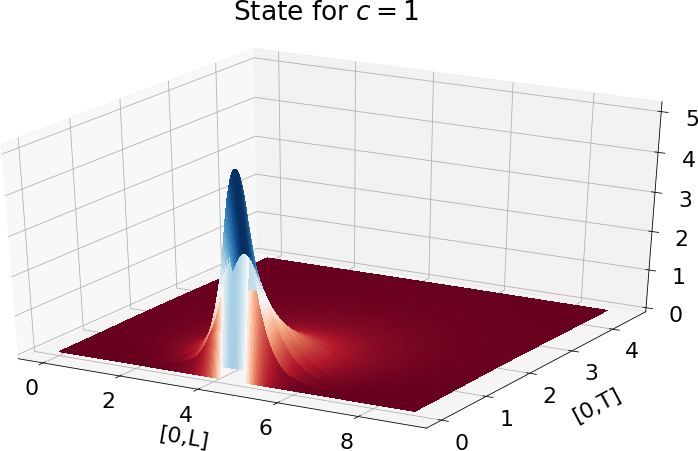}
	\caption{Influence of perturbations in the uncontrolled case.}
	\label{fig:parab_sim}
\end{figure}

In Figure~\ref{fig:parab_sim}, we depict the results of a simulation with a perturbation concentrated at the beginning of the temporal domain and in the middle of the spatial domain, i.e.,
\begin{align*}
f(t,\omega) = \begin{cases}
10 \qquad & \omega\in L/2\ \wedge t\in [0,2]\\
0 \qquad & \mathrm{otherwise}.
\end{cases} 
\end{align*}
We observe that, as in the stationary case, the influence of the perturbation has a very different effect in the three uncontrolled equations. For $c=0$, we obtain the heat equation which converges to a non-zero steady state due to the homogeneous Neumann boundary conditions and the resulting one-dimensional kernel spanned by the constant function. For $c=-1$, the zero eigenvalue is shifted by one to the right such that the equation is exponentially unstable. Contrary, for $c=1$, the equation is exponentially stable.

\begin{figure}[H]
	\centering
	\includegraphics[width=.32\linewidth]{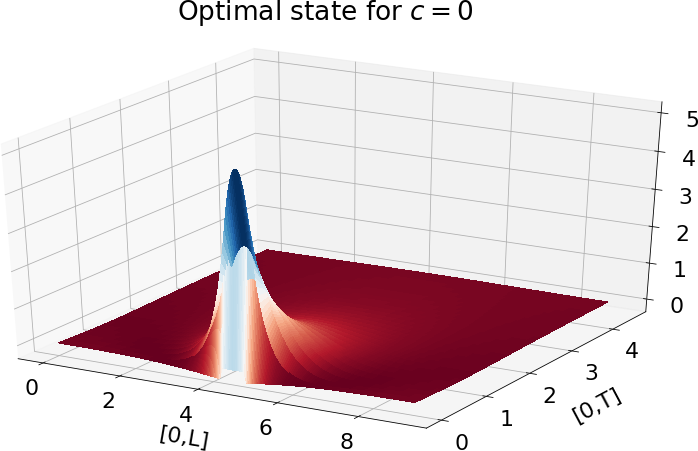}
	\includegraphics[width=.32\linewidth]{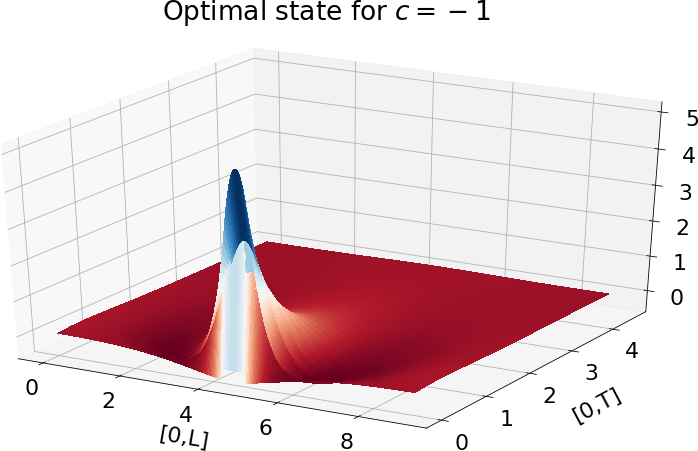}
	\includegraphics[width=.32\linewidth]{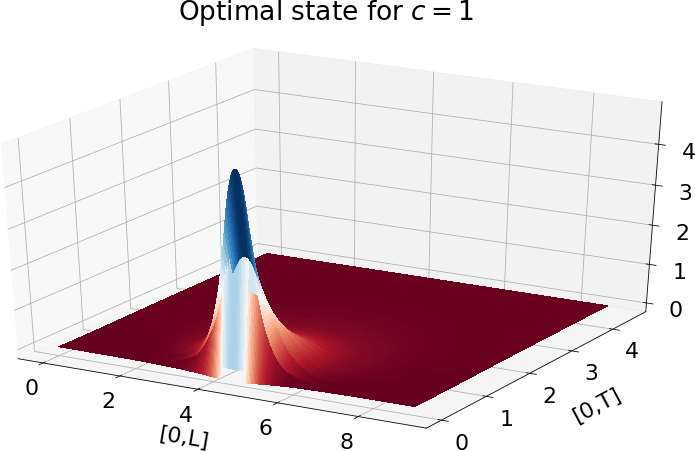}
	\caption{Influence of perturbations on the optimal state.}
	\label{fig:parab_oc}
\end{figure}

\noindent Second we present the optimal states for optimal control with the cost functional
\begin{align*}
\min_{u\in L^2(0,T;L^2(\Omega))} \frac12\left(  \|y\|^2_{L^2(0,T;L^2(\Omega))}+ \|u\|^2_{L^2(0,T;L^2(\Omega))}  \right).
\end{align*}
As a particular case of Example~\ref{ex:laplace} and as can be seen directly by choosing the feedback $K_B y = K_C^*y = cy - \varepsilon y$ for any $\varepsilon>0$, this system satisfies the assumptions of Theorem~\ref{thm:opnorm_p} and hence, we obtain an exponential decay of perturbations in view of Theorem~\ref{thm:scaling_p}. This can be clearly, observed in Figure~\ref{fig:parab_oc}, i.e., all three equations converge to the origin. The convergence behavior in time is implied by the turnpike property and can be characterized by means of Riccati equations~\cite{TrelZuaz15}.
A similar characterization of the decay in space goes beyond the scope of this work but is subject to future research.

\section{Conclusion and future work}\label{sec:conclusions}
\noindent In this work we provided sufficient conditions to deduce exponential decay of perturbations in space for elliptic and parabolic optimal control problems. 
To this end, we first proved that the solution operator of the optimality system is uniformly bounded in the domain size provided that certain stabilizability and detectability conditions hold. 
	Then, we leveraged this key property to show that optimality implies exponential localization, meaning that the influence of perturbations on optimal solutions decays exponentially in space. 
To this end, we combined a novel scaling result with a Neumann argument.
Further, we discussed various examples showing the straightforward verifiability of the assumptions for different 
settings, including boundary control, distributed control, and coupled problems.

Future work might consider the extension to hyperbolic equations and networks thereof, showing that errors also only act local in the network. 
\bibliographystyle{abbrv}
\bibliography{references}

\begin{thebibliography}{10}

\bibitem{AlnaBlech15}
M.~Aln{\ae}s, J.~Blechta, J.~Hake, A.~Johansson, B.~Kehlet, A.~Logg,
  C.~Richardson, J.~Ring, M.~E. Rognes, and G.~N. Wells.
\newblock The {FEniCS} project version 1.5.
\newblock {\em Archive of Numerical Software}, 3(100), 2015.

\bibitem{BartHein06}
R.~A. Bartlett, M.~Heinkenschloss, D.~Ridzal, and B.~G. van Bloemen~Waanders.
\newblock Domain decomposition methods for advection dominated linear-quadratic
  elliptic optimal control problems.
\newblock {\em Computer methods in applied mechanics and engineering},
  195(44-47):6428--6447, 2006.

\bibitem{Bebe03}
M.~Bebendorf.
\newblock A note on the {P}oincar{\'e} inequality for convex domains.
\newblock {\em Zeitschrift f{\"u}r Analysis und ihre Anwendungen},
  22(4):751--756, 2003.

\bibitem{Brezi11}
H.~Br{\'e}zis.
\newblock {\em Functional analysis, Sobolev spaces and partial differential
  equations}, volume~2.
\newblock Springer, 2011.

\bibitem{DammGrun14}
T.~Damm, L.~Grüne, M.~Stieler, and K.~Worthmann.
\newblock An exponential turnpike theorem for dissipative discrete time optimal
  control problems.
\newblock {\em SIAM Journal on Control and Optimization}, 52(3):1935--1957,
  2014.

\bibitem{EsteMele11}
S.~Esterhazy and J.~M. Melenk.
\newblock On stability of discretizations of the {H}elmholtz equation.
\newblock In {\em Numerical analysis of multiscale problems}, pages 285--324.
  Springer, 2011.

\bibitem{EsteGesh20}
C.~Esteve, B.~Geshkovski, D.~Pighin, and E.~Zuazua.
\newblock Large-time asymptotics in deep learning.
\newblock {\em preprint arXiv:2008.02491}, 2020.

\bibitem{evans2022partial}
L.~C. Evans.
\newblock {\em Partial differential equations}, volume~19.
\newblock American Mathematical Society, 2022.

\bibitem{FaulFlas22}
T.~Faulwasser, K.~Fla{\ss}kamp, S.~Ober-Bl{\"o}baum, M.~Schaller, and
  K.~Worthmann.
\newblock Manifold turnpikes, trims, and symmetries.
\newblock {\em Mathematics of Control, Signals, and Systems}, 34(4):759--788,
  2022.

\bibitem{FaulGrun22}
T.~Faulwasser and L.~Gr{\"u}ne.
\newblock Turnpike properties in optimal control: An overview of discrete-time
  and continuous-time results.
\newblock {\em Handbook of numerical analysis}, 23:367--400, 2022.

\bibitem{FaulHemp21}
T.~Faulwasser, A.-J. Hempel, and S.~Streif.
\newblock On the turnpike to design of deep neural nets: Explicit depth bounds.
\newblock {\em arXiv preprint arXiv:2101.03000}, 2021.

\bibitem{FaulMurr20}
T.~Faulwasser and A.~Murray.
\newblock Turnpike properties in discrete-time mixed-integer optimal control.
\newblock {\em IEEE Control Systems Letters}, 4(3):704--709, 2020.

\bibitem{GeshZuaz22}
B.~Geshkovski and E.~Zuazua.
\newblock Turnpike in optimal control of {PDEs}, {ResNets}, and beyond.
\newblock {\em Acta Numerica}, 31:135--263, 2022.

\bibitem{Grun13}
L.~Gr{\"u}ne.
\newblock Economic receding horizon control without terminal constraints.
\newblock {\em Automatica}, 49(3):725--734, 2013.

\bibitem{Grun16}
L.~Gr{\"u}ne.
\newblock Approximation properties of receding horizon optimal control.
\newblock {\em Jahresbericht der Deutschen Mathematiker-Vereinigung},
  118:3--37, 2016.

\bibitem{GrunScha20}
L.~Gr{\"u}ne, M.~Schaller, and A.~Schiela.
\newblock Exponential sensitivity and turnpike analysis for linear quadratic
  optimal control of general evolution equations.
\newblock {\em Journal of Differential Equations}, 268(12):7311--7341, 2020.

\bibitem{GrunScha21}
L.~Gr{\"u}ne, M.~Schaller, and A.~Schiela.
\newblock Abstract nonlinear sensitivity and turnpike analysis and an
  application to semilinear parabolic {PDE}s.
\newblock {\em ESAIM: Control, Optimisation and Calculus of Variations}, 27:56,
  2021.

\bibitem{GrunScha19}
L.~Grüne, M.~Schaller, and A.~Schiela.
\newblock Sensitivity analysis of optimal control for a class of parabolic
  {PDE}s motivated by model predictive control.
\newblock {\em SIAM Journal on Control and Optimization}, 57(4):2753--2774,
  2019.

\bibitem{GrunScha22}
L.~Grüne, M.~Schaller, and A.~Schiela.
\newblock Efficient model predictive control for parabolic {PDE}s with goal
  oriented error estimation.
\newblock {\em SIAM Journal on Scientific Computing}, 44(1):A471--A500, 2022.

\bibitem{GugaTrel16}
M.~Gugat, E.~Tr{\'e}lat, and E.~Zuazua.
\newblock Optimal neumann control for the 1d wave equation: finite horizon,
  infinite horizon, boundary tracking terms and the turnpike property.
\newblock {\em Systems \& Control Letters}, 90:61--70, 2016.

\bibitem{GuglLi24}
R.~Guglielmi and Z.~Li.
\newblock Necessary conditions for turnpike property for generalized
  linear--quadratic problems.
\newblock {\em Mathematics of Control, Signals, and Systems}, pages 1--31,
  2024.

\bibitem{HaucPete22}
M.~Hauck and D.~Peterseim.
\newblock Multi-resolution localized orthogonal decomposition for {H}elmholtz
  problems.
\newblock {\em Multiscale Modeling \& Simulation}, 20(2):657--684, 2022.

\bibitem{LagnLeug04}
J.~E. Lagnese and G.~Leugering.
\newblock {\em Domain in Decomposition Methods in Optimal Control of Partial
  Differential Equations}.
\newblock Number 148. Springer Science \& Business Media, 2004.

\bibitem{NaShin22}
S.~Na, S.~Shin, M.~Anitescu, and V.~M. Zavala.
\newblock On the convergence of overlapping {S}chwarz decomposition for
  nonlinear optimal control.
\newblock {\em IEEE Transactions on Automatic Control}, 67(11):5996--6011,
  2022.

\bibitem{Scha21}
M.~Schaller.
\newblock {\em Sensitivity Analysis and Goal Oriented Error Estimation for
  Model Predictive Control}.
\newblock PhD thesis, University of Bayreuth, 2021.

\bibitem{Schi13}
A.~Schiela.
\newblock A concise proof for existence and uniqueness of solutions of linear
  parabolic {PDE}s in the context of optimal control.
\newblock {\em Systems \& Control Letters}, 62(10):895--901, 2013.

\bibitem{Shin21}
S.~Shin.
\newblock {\em Graph-Structured Nonlinear Programming: Properties and
  Algorithms}.
\newblock PhD thesis, The University of Wisconsin-Madison, 2021.

\bibitem{ShinAnit22}
S.~Shin, M.~Anitescu, and V.~M. Zavala.
\newblock Exponential decay of sensitivity in graph-structured nonlinear
  programs.
\newblock {\em SIAM Journal on Optimization}, 32(2):1156--1183, 2022.

\bibitem{ShinFaul19}
S.~Shin, T.~Faulwasser, M.~Zanon, and V.~M. Zavala.
\newblock A parallel decomposition scheme for solving long-horizon optimal
  control problems.
\newblock In {\em 2019 IEEE 58th Conference on Decision and Control (CDC)},
  pages 5264--5271. IEEE, 2019.

\bibitem{ShinLin23}
S.~Shin, Y.~Lin, G.~Qu, A.~Wierman, and M.~Anitescu.
\newblock Near-optimal distributed linear-quadratic regulator for networked
  systems.
\newblock {\em SIAM Journal on Control and Optimization}, 61(3):1113--1135,
  2023.

\bibitem{ShinZava21}
S.~Shin and V.~M. Zavala.
\newblock Diffusing-horizon model predictive control.
\newblock {\em IEEE Transactions on Automatic Control}, 2021.

\bibitem{SperSalu23}
M.~Sperl, L.~Saluzzi, L.~Gr{\"u}ne, and D.~Kalise.
\newblock Separable approximations of optimal value functions under a decaying
  sensitivity assumption.
\newblock In {\em Proceedings of the 62nd IEEE Conference on Decision and
  Control (CDC)}, pages 259--264, 2023.

\bibitem{StakHols11}
I.~Stakgold and M.~J. Holst.
\newblock {\em Green's functions and boundary value problems}.
\newblock John Wiley \& Sons, 2011.

\bibitem{TrelZuaz15}
E.~Tr{\'e}lat and E.~Zuazua.
\newblock The turnpike property in finite-dimensional nonlinear optimal
  control.
\newblock {\em Journal of Differential Equations}, 258(1):81--114, 2015.

\bibitem{VirtGomm20}
P.~Virtanen, R.~Gommers, T.~E. Oliphant, M.~Haberland, T.~Reddy, D.~Cournapeau,
  E.~Burovski, P.~Peterson, W.~Weckesser, J.~Bright, S.~J. {van der Walt},
  M.~Brett, J.~Wilson, K.~J. Millman, N.~Mayorov, A.~R.~J. Nelson, E.~Jones,
  R.~Kern, E.~Larson, C.~J. Carey, {\.I}.~Polat, Y.~Feng, E.~W. Moore,
  J.~{VanderPlas}, D.~Laxalde, J.~Perktold, R.~Cimrman, I.~Henriksen, E.~A.
  Quintero, C.~R. Harris, A.~M. Archibald, A.~H. Ribeiro, F.~Pedregosa, P.~{van
  Mulbregt}, and {SciPy 1.0 Contributors}.
\newblock {{SciPy} 1.0: Fundamental Algorithms for Scientific Computing in
  Python}.
\newblock {\em Nature Methods}, 17:261--272, 2020.

\bibitem{Wlok87}
J.~Wloka.
\newblock {\em Partial differential equations}.
\newblock Cambridge University Press, 1987.

\bibitem{yang2019survey}
T.~Yang, X.~Yi, J.~Wu, Y.~Yuan, D.~Wu, Z.~Meng, Y.~Hong, H.~Wang, Z.~Lin, and
  K.~H. Johansson.
\newblock A survey of distributed optimization.
\newblock {\em Annual Reviews in Control}, 47:278--305, 2019.

\bibitem{Zeid13}
E.~Zeidler.
\newblock {\em Nonlinear functional analysis and its applications: II/B:
  Nonlinear monotone operators}.
\newblock Springer Science \& Business Media, 2013.

\end{thebibliography}
\end{document}